\documentclass{gtart}

\def\ifplaintex{\expandafter\ifx\csname documentclass\endcsname\relax}

\def\gtp{{\mathsurround=0pt\it $\cal G\mskip-2mu$eometry \&\ 
$\cal T\!\!$opology $\cal P\!$ublications}}  

\def\Addressesr{\bigskip
{\small \parskip 0pt \leftskip 0pt \rightskip 0pt plus 1fil \def\\{\par}
\sl\theaddress\par
\medskip
\rm Email:\stdspace\tt\theemail\hfill\rm Received:\qua\receiveddate \par}}

\def\recd{{\small Received:\qua\receiveddate\ifx\reviseddate\relax
\else\qquad Revised:\qua\reviseddate\fi\par}} 

\def\lognumber#1{\def\thelognumber{#1}}
\def\volumenumber#1{\def\thevolumenumber{#1}}
\def\volumeyear#1{\def\thevolumeyear{#1}}
\def\papernumber#1{\def\thepapernumber{#1}}
\def\pagenumbers#1#2{\def\startpage{#1}\def\finishpage{#2}}
\def\published#1{\def\publishdate{#1}}

\def\received#1{\def\receiveddate{#1}}

\def\accepted#1{\def\accepteddate{#1}}

\def\asciiauthors#1{\def\theasciiauthors{#1}}

\def\coverauthors#1{\def\thecoverauthors{#1}}
\long\def\asciiabstract#1{\long\def\theasciiabstract{#1}}

\let\\\par\let\thelognumber\relax\let\thevolumenumber\relax
\let\thepapernumber\relax\let\thevolumeyear\relax\let\startpage\relax
\let\finishpage\relax\let\publishdate\relax\let\receiveddate\relax
\let\reviseddate\relax\let\accepteddate\relax\let\theasciititle\relax
\let\theasciiauthors\relax
\let\theasciiabstract\relax

\let\thecoverauthors\relax\let\theasciiemail\relax

\ifplaintex
\font\logobig=cmssbx10 scaled 3836
\font\logomed=cmssbx10 scaled 2557
\else
\font\logobig=cmssbx10 scaled 4200
\font\logomed=cmssbx10 scaled 2800
\fi

\long\def\makeagttitle{   
\count0=\startpage
\agt\hfill      
\hbox to 45truept{\vbox to 0pt{\vglue -13truept{\logomed A\kern -.37em{\logobig 
T}\kern -.38em G}\vss}\hss}
\break
{\small Volume \thevolumenumber\ (\thevolumeyear)
\startpage--\finishpage\nl
Published: \publishdate}

\vglue .25truein

{\parskip=0pt\leftskip 0pt plus
1fil\def\\{\par\smallskip}{\Large\bf\thetitle}\par\medskip} \vglue
0.05truein

{\parskip=0pt\leftskip 0pt plus 1fil\def\\{\par}{\sc\theauthors}
\par\medskip}%
 
\vglue 0.03truein 

{\small\leftskip 25truept\rightskip 25truept{\bf Abstract}\stdspace\theabstract

{\bf AMS Classification}\stdspace\theprimaryclass
\ifx\thesecondaryclass\relax\else; \thesecondaryclass\fi\par
{\bf Keywords}\stdspace \thekeywords\par}\vglue 7truept

}   

\ifplaintex
\font\phead=cmsl9 scaled 950
\font\pnum=cmbx10 scaled 913
\font\pfoot=cmsl9 scaled 950
\headline{\vbox to 0pt{\vskip -4.5mm\line{\small\phead\ifnum
\count0=\startpage ISSN 1472-2739 (on-line) 1472-2747 (printed)
\hfill {\pnum\folio}\else\ifodd\count0\def\\{ }%
\ifx\theshorttitle\relax\thetitle\else\theshorttitle\fi\hfill{\pnum\folio}
\else\def\\{ and }{\pnum\folio}\hfill\ifx\theshortauthors\relax\theauthors
\else\theshortauthors\fi\fi\fi}\vss}}
\footline{\vbox to 0pt{\vglue 0mm\line{\small\pfoot\ifnum\count0=\startpage
\copyright\ \gtp\hfill\else
\agt, Volume \thevolumenumber\ (\thevolumeyear)\hfill\fi}\vss}}
\else
\font=cmsl9 scaled 1050
\font\lnum=cmbx10 
\font=cmsl9 scaled 1050
\makeatletter
\def\@oddhead{{\small\lhead\ifnum\count0=\startpage ISSN 1472-2739 
(on-line) 1472-2747 (printed)\hfill {\lnum\number\count0}\else\ifodd\count0
\def\\{ }\ifx\theshorttitle\relax \thetitle \else\theshorttitle\fi\hfill
{\lnum\number\count0}\else\def\\{ and }{\lnum\number\count0}
\hfill\ifx\theshortauthors\relax 
\theauthors\else\theshortauthors\fi\fi\fi}}\def\@evenhead{\@oddhead}
\def\@oddfoot{\small\lfoot\ifnum\count0=\startpage\copyright\ \gtp\hfill\else
\agt, Volume \thevolumenumber\ (\thevolumeyear)\hfill\fi}
\def\@evenfoot{\@oddfoot}
\makeatother
\fi
\let\maketitlepage\makeagttitle
\let\makeshorttitle\maketitlepage
\let\maketitle\maketitlepage

\newwrite\gtoutfile
\long\gdef\makeheadfile{  
{\def\\{, }\def\s{ }
\immediate\openout\gtoutfile head.xxx
\immediate\write\gtoutfile{To: math@arxiv.org}
\immediate\write\gtoutfile{Subject: put OR rep NNNNN:ppppp}
\immediate\write\gtoutfile{--text follows this line--}
\immediate\write\gtoutfile{Proxy-for: \ifx\theasciiauthors\relax
\theauthors\else\theasciiauthors\fi\s<\ifx\theasciiemail\relax\theemail\else\theasciiemail\fi>}
\immediate\write\gtoutfile{\noexpand\\}
\immediate\write\gtoutfile{Authors: \ifx\theasciiauthors\relax
\theauthors\else\theasciiauthors\fi}
{\def\\{ }\immediate\write\gtoutfile{Title: \ifx\theasciititle\relax
\thetitle\else\theasciititle\fi}}
\immediate\write\gtoutfile{Subj-class: GT or SG, GR etc}
\immediate\write\gtoutfile{MSC-class: \theprimaryclass\ifx\thesecondaryclass\relax\else, \thesecondaryclass\fi}
\immediate\write\gtoutfile{Journal-ref: Algebr. Geom. Topol. \thevolumenumber\s
(\thevolumeyear) \startpage-\finishpage}
\immediate\write\gtoutfile{Comments: Published by Algebraic and
Geometric Topology at}
\immediate\write\gtoutfile{\s\s\s  http://www.maths.warwick.ac.uk/agt/AGTVol\thevolumenumber/agt-\thevolumenumber-\thepapernumber.abs.html}
\immediate\write\gtoutfile{\noexpand\\}
\immediate\write\gtoutfile{}
\ifx\theasciiabstract\relax
\immediate\write\gtoutfile{\theabstract}\else
\immediate\write\gtoutfile{\theasciiabstract}\fi
\immediate\write\gtoutfile{}
\immediate\write\gtoutfile{\noexpand\\}
\immediate\write\gtoutfile{}
\immediate\closeout\gtoutfile}}  

\def\maketitlepage{\makeagttitle\makeheadfile}
\let\makeshorttitle\maketitlepage
\let\maketitle\maketitlepage

\def\ifplaintex{\expandafter\ifx\csname documentclass\endcsname\relax}

\def\gtp{{\mathsurround=0pt\it $\cal G\mskip-2mu$eometry \&\ 
$\cal T\!\!$opology $\cal P\!$ublications}}  

\def\Addressesr{\bigskip
{\small \parskip 0pt \leftskip 0pt \rightskip 0pt plus 1fil \def\\{\par}
\sl\theaddress\par
\medskip
\rm Email:\stdspace\tt\theemail\hfill\rm Received:\qua\receiveddate \par}}

\def\recd{{\small Received:\qua\receiveddate\ifx\reviseddate\relax
\else\qquad Revised:\qua\reviseddate\fi\par}} 

\def\lognumber#1{\def\thelognumber{#1}}
\def\volumenumber#1{\def\thevolumenumber{#1}}
\def\volumeyear#1{\def\thevolumeyear{#1}}
\def\papernumber#1{\def\thepapernumber{#1}}
\def\pagenumbers#1#2{\def\startpage{#1}\def\finishpage{#2}}
\def\published#1{\def\publishdate{#1}}

\def\received#1{\def\receiveddate{#1}}

\def\accepted#1{\def\accepteddate{#1}}

\def\asciiauthors#1{\def\theasciiauthors{#1}}

\def\coverauthors#1{\def\thecoverauthors{#1}}
\long\def\asciiabstract#1{\long\def\theasciiabstract{#1}}

\let\\\par\let\thelognumber\relax\let\thevolumenumber\relax
\let\thepapernumber\relax\let\thevolumeyear\relax\let\startpage\relax
\let\finishpage\relax\let\publishdate\relax\let\receiveddate\relax
\let\reviseddate\relax\let\accepteddate\relax\let\theasciititle\relax
\let\theasciiauthors\relax
\let\theasciiabstract\relax

\let\thecoverauthors\relax\let\theasciiemail\relax

\ifplaintex
\font\logobig=cmssbx10 scaled 3836
\font\logomed=cmssbx10 scaled 2557
\else
\font\logobig=cmssbx10 scaled 4200
\font\logomed=cmssbx10 scaled 2800
\fi

\long\def\makeagttitle{   
\count0=\startpage
\agt\hfill      
\hbox to 45truept{\vbox to 0pt{\vglue -13truept{\logomed A\kern -.37em{\logobig 
T}\kern -.38em G}\vss}\hss}
\break
{\small Volume \thevolumenumber\ (\thevolumeyear)
\startpage--\finishpage\nl
Published: \publishdate}

\vglue .25truein

{\parskip=0pt\leftskip 0pt plus
1fil\def\\{\par\smallskip}{\Large\bf\thetitle}\par\medskip} \vglue
0.05truein

{\parskip=0pt\leftskip 0pt plus 1fil\def\\{\par}{\sc\theauthors}
\par\medskip}%
 
\vglue 0.03truein 

{\small\leftskip 25truept\rightskip 25truept{\bf Abstract}\stdspace\theabstract

{\bf AMS Classification}\stdspace\theprimaryclass
\ifx\thesecondaryclass\relax\else; \thesecondaryclass\fi\par
{\bf Keywords}\stdspace \thekeywords\par}\vglue 7truept

}   

\ifplaintex
\font\phead=cmsl9 scaled 950
\font\pnum=cmbx10 scaled 913
\font\pfoot=cmsl9 scaled 950
\headline{\vbox to 0pt{\vskip -4.5mm\line{\small\phead\ifnum
\count0=\startpage ISSN 1472-2739 (on-line) 1472-2747 (printed)
\hfill {\pnum\folio}\else\ifodd\count0\def\\{ }%
\ifx\theshorttitle\relax\thetitle\else\theshorttitle\fi\hfill{\pnum\folio}
\else\def\\{ and }{\pnum\folio}\hfill\ifx\theshortauthors\relax\theauthors
\else\theshortauthors\fi\fi\fi}\vss}}
\footline{\vbox to 0pt{\vglue 0mm\line{\small\pfoot\ifnum\count0=\startpage
\copyright\ \gtp\hfill\else
\agt, Volume \thevolumenumber\ (\thevolumeyear)\hfill\fi}\vss}}
\else
\font=cmsl9 scaled 1050
\font\lnum=cmbx10 
\font=cmsl9 scaled 1050
\makeatletter
\def\@oddhead{{\small\lhead\ifnum\count0=\startpage ISSN 1472-2739 
(on-line) 1472-2747 (printed)\hfill {\lnum\number\count0}\else\ifodd\count0
\def\\{ }\ifx\theshorttitle\relax \thetitle \else\theshorttitle\fi\hfill
{\lnum\number\count0}\else\def\\{ and }{\lnum\number\count0}
\hfill\ifx\theshortauthors\relax 
\theauthors\else\theshortauthors\fi\fi\fi}}\def\@evenhead{\@oddhead}
\def\@oddfoot{\small\lfoot\ifnum\count0=\startpage\copyright\ \gtp\hfill\else
\agt, Volume \thevolumenumber\ (\thevolumeyear)\hfill\fi}
\def\@evenfoot{\@oddfoot}
\makeatother
\fi
\let\maketitlepage\makeagttitle
\let\makeshorttitle\maketitlepage
\let\maketitle\maketitlepage

\newwrite\gtoutfile
\long\gdef\makeheadfile{  
{\def\\{, }\def\s{ }
\immediate\openout\gtoutfile head.xxx
\immediate\write\gtoutfile{To: math@arxiv.org}
\immediate\write\gtoutfile{Subject: put OR rep NNNNN:ppppp}
\immediate\write\gtoutfile{--text follows this line--}
\immediate\write\gtoutfile{Proxy-for: \ifx\theasciiauthors\relax
\theauthors\else\theasciiauthors\fi\s<\ifx\theasciiemail\relax\theemail\else\theasciiemail\fi>}
\immediate\write\gtoutfile{\noexpand\\}
\immediate\write\gtoutfile{Authors: \ifx\theasciiauthors\relax
\theauthors\else\theasciiauthors\fi}
{\def\\{ }\immediate\write\gtoutfile{Title: \ifx\theasciititle\relax
\thetitle\else\theasciititle\fi}}
\immediate\write\gtoutfile{Subj-class: GT or SG, GR etc}
\immediate\write\gtoutfile{MSC-class: \theprimaryclass\ifx\thesecondaryclass\relax\else, \thesecondaryclass\fi}
\immediate\write\gtoutfile{Journal-ref: Algebr. Geom. Topol. \thevolumenumber\s
(\thevolumeyear) \startpage-\finishpage}
\immediate\write\gtoutfile{Comments: Published by Algebraic and
Geometric Topology at}
\immediate\write\gtoutfile{\s\s\s  http://www.maths.warwick.ac.uk/agt/AGTVol\thevolumenumber/agt-\thevolumenumber-\thepapernumber.abs.html}
\immediate\write\gtoutfile{\noexpand\\}
\immediate\write\gtoutfile{}
\ifx\theasciiabstract\relax
\immediate\write\gtoutfile{\theabstract}\else
\immediate\write\gtoutfile{\theasciiabstract}\fi
\immediate\write\gtoutfile{}
\immediate\write\gtoutfile{\noexpand\\}
\immediate\write\gtoutfile{}
\immediate\closeout\gtoutfile}}  

\def\maketitlepage{\makeagttitle\makeheadfile}
\let\makeshorttitle\maketitlepage
\let\maketitle\maketitlepage

\def\ifplaintex{\expandafter\ifx\csname documentclass\endcsname\relax}

\def\gtp{{\mathsurround=0pt\it $\cal G\mskip-2mu$eometry \&\ 
$\cal T\!\!$opology $\cal P\!$ublications}}  

\def\Addressesr{\bigskip
{\small \parskip 0pt \leftskip 0pt \rightskip 0pt plus 1fil \def\\{\par}
\sl\theaddress\par
\medskip
\rm Email:\stdspace\tt\theemail\hfill\rm Received:\qua\receiveddate \par}}

\def\recd{{\small Received:\qua\receiveddate\ifx\reviseddate\relax
\else\qquad Revised:\qua\reviseddate\fi\par}} 

\def\lognumber#1{\def\thelognumber{#1}}
\def\volumenumber#1{\def\thevolumenumber{#1}}
\def\volumeyear#1{\def\thevolumeyear{#1}}
\def\papernumber#1{\def\thepapernumber{#1}}
\def\pagenumbers#1#2{\def\startpage{#1}\def\finishpage{#2}}
\def\published#1{\def\publishdate{#1}}

\def\received#1{\def\receiveddate{#1}}

\def\accepted#1{\def\accepteddate{#1}}

\def\asciiauthors#1{\def\theasciiauthors{#1}}

\def\coverauthors#1{\def\thecoverauthors{#1}}
\long\def\asciiabstract#1{\long\def\theasciiabstract{#1}}

\let\\\par\let\thelognumber\relax\let\thevolumenumber\relax
\let\thepapernumber\relax\let\thevolumeyear\relax\let\startpage\relax
\let\finishpage\relax\let\publishdate\relax\let\receiveddate\relax
\let\reviseddate\relax\let\accepteddate\relax\let\theasciititle\relax
\let\theasciiauthors\relax
\let\theasciiabstract\relax

\let\thecoverauthors\relax\let\theasciiemail\relax

\ifplaintex
\font\logobig=cmssbx10 scaled 3836
\font\logomed=cmssbx10 scaled 2557
\else
\font\logobig=cmssbx10 scaled 4200
\font\logomed=cmssbx10 scaled 2800
\fi

\long\def\makeagttitle{   
\count0=\startpage
\agt\hfill      
\hbox to 45truept{\vbox to 0pt{\vglue -13truept{\logomed A\kern -.37em{\logobig 
T}\kern -.38em G}\vss}\hss}
\break
{\small Volume \thevolumenumber\ (\thevolumeyear)
\startpage--\finishpage\nl
Published: \publishdate}

\vglue .25truein

{\parskip=0pt\leftskip 0pt plus
1fil\def\\{\par\smallskip}{\Large\bf\thetitle}\par\medskip} \vglue
0.05truein

{\parskip=0pt\leftskip 0pt plus 1fil\def\\{\par}{\sc\theauthors}
\par\medskip}%
 
\vglue 0.03truein 

{\small\leftskip 25truept\rightskip 25truept{\bf Abstract}\stdspace\theabstract

{\bf AMS Classification}\stdspace\theprimaryclass
\ifx\thesecondaryclass\relax\else; \thesecondaryclass\fi\par
{\bf Keywords}\stdspace \thekeywords\par}\vglue 7truept

}   

\ifplaintex
\font\phead=cmsl9 scaled 950
\font\pnum=cmbx10 scaled 913
\font\pfoot=cmsl9 scaled 950
\headline{\vbox to 0pt{\vskip -4.5mm\line{\small\phead\ifnum
\count0=\startpage ISSN 1472-2739 (on-line) 1472-2747 (printed)
\hfill {\pnum\folio}\else\ifodd\count0\def\\{ }%
\ifx\theshorttitle\relax\thetitle\else\theshorttitle\fi\hfill{\pnum\folio}
\else\def\\{ and }{\pnum\folio}\hfill\ifx\theshortauthors\relax\theauthors
\else\theshortauthors\fi\fi\fi}\vss}}
\footline{\vbox to 0pt{\vglue 0mm\line{\small\pfoot\ifnum\count0=\startpage
\copyright\ \gtp\hfill\else
\agt, Volume \thevolumenumber\ (\thevolumeyear)\hfill\fi}\vss}}
\else
\font=cmsl9 scaled 1050
\font\lnum=cmbx10 
\font=cmsl9 scaled 1050
\makeatletter
\def\@oddhead{{\small\lhead\ifnum\count0=\startpage ISSN 1472-2739 
(on-line) 1472-2747 (printed)\hfill {\lnum\number\count0}\else\ifodd\count0
\def\\{ }\ifx\theshorttitle\relax \thetitle \else\theshorttitle\fi\hfill
{\lnum\number\count0}\else\def\\{ and }{\lnum\number\count0}
\hfill\ifx\theshortauthors\relax 
\theauthors\else\theshortauthors\fi\fi\fi}}\def\@evenhead{\@oddhead}
\def\@oddfoot{\small\lfoot\ifnum\count0=\startpage\copyright\ \gtp\hfill\else
\agt, Volume \thevolumenumber\ (\thevolumeyear)\hfill\fi}
\def\@evenfoot{\@oddfoot}
\makeatother
\fi
\let\maketitlepage\makeagttitle
\let\makeshorttitle\maketitlepage
\let\maketitle\maketitlepage

\newwrite\gtoutfile
\long\gdef\makeheadfile{  
{\def\\{, }\def\s{ }
\immediate\openout\gtoutfile head.xxx
\immediate\write\gtoutfile{To: math@arxiv.org}
\immediate\write\gtoutfile{Subject: put OR rep NNNNN:ppppp}
\immediate\write\gtoutfile{--text follows this line--}
\immediate\write\gtoutfile{Proxy-for: \ifx\theasciiauthors\relax
\theauthors\else\theasciiauthors\fi\s<\ifx\theasciiemail\relax\theemail\else\theasciiemail\fi>}
\immediate\write\gtoutfile{\noexpand\\}
\immediate\write\gtoutfile{Authors: \ifx\theasciiauthors\relax
\theauthors\else\theasciiauthors\fi}
{\def\\{ }\immediate\write\gtoutfile{Title: \ifx\theasciititle\relax
\thetitle\else\theasciititle\fi}}
\immediate\write\gtoutfile{Subj-class: GT or SG, GR etc}
\immediate\write\gtoutfile{MSC-class: \theprimaryclass\ifx\thesecondaryclass\relax\else, \thesecondaryclass\fi}
\immediate\write\gtoutfile{Journal-ref: Algebr. Geom. Topol. \thevolumenumber\s
(\thevolumeyear) \startpage-\finishpage}
\immediate\write\gtoutfile{Comments: Published by Algebraic and
Geometric Topology at}
\immediate\write\gtoutfile{\s\s\s  http://www.maths.warwick.ac.uk/agt/AGTVol\thevolumenumber/agt-\thevolumenumber-\thepapernumber.abs.html}
\immediate\write\gtoutfile{\noexpand\\}
\immediate\write\gtoutfile{}
\ifx\theasciiabstract\relax
\immediate\write\gtoutfile{\theabstract}\else
\immediate\write\gtoutfile{\theasciiabstract}\fi
\immediate\write\gtoutfile{}
\immediate\write\gtoutfile{\noexpand\\}
\immediate\write\gtoutfile{}
\immediate\closeout\gtoutfile}}  

\def\maketitlepage{\makeagttitle\makeheadfile}
\let\makeshorttitle\maketitlepage
\let\maketitle\maketitlepage

\def\ifplaintex{\expandafter\ifx\csname documentclass\endcsname\relax}

\def\gtp{{\mathsurround=0pt\it $\cal G\mskip-2mu$eometry \&\ 
$\cal T\!\!$opology $\cal P\!$ublications}}  

\def\Addressesr{\bigskip
{\small \parskip 0pt \leftskip 0pt \rightskip 0pt plus 1fil \def\\{\par}
\sl\theaddress\par
\medskip
\rm Email:\stdspace\tt\theemail\hfill\rm Received:\qua\receiveddate \par}}

\def\recd{{\small Received:\qua\receiveddate\ifx\reviseddate\relax
\else\qquad Revised:\qua\reviseddate\fi\par}} 

\def\lognumber#1{\def\thelognumber{#1}}
\def\volumenumber#1{\def\thevolumenumber{#1}}
\def\volumeyear#1{\def\thevolumeyear{#1}}
\def\papernumber#1{\def\thepapernumber{#1}}
\def\pagenumbers#1#2{\def\startpage{#1}\def\finishpage{#2}}
\def\published#1{\def\publishdate{#1}}

\def\received#1{\def\receiveddate{#1}}

\def\accepted#1{\def\accepteddate{#1}}

\def\asciiauthors#1{\def\theasciiauthors{#1}}

\def\coverauthors#1{\def\thecoverauthors{#1}}
\long\def\asciiabstract#1{\long\def\theasciiabstract{#1}}

\let\\\par\let\thelognumber\relax\let\thevolumenumber\relax
\let\thepapernumber\relax\let\thevolumeyear\relax\let\startpage\relax
\let\finishpage\relax\let\publishdate\relax\let\receiveddate\relax
\let\reviseddate\relax\let\accepteddate\relax\let\theasciititle\relax
\let\theasciiauthors\relax
\let\theasciiabstract\relax

\let\thecoverauthors\relax\let\theasciiemail\relax

\ifplaintex
\font\logobig=cmssbx10 scaled 3836
\font\logomed=cmssbx10 scaled 2557
\else
\font\logobig=cmssbx10 scaled 4200
\font\logomed=cmssbx10 scaled 2800
\fi

\long\def\makeagttitle{   
\count0=\startpage
\agt\hfill      
\hbox to 45truept{\vbox to 0pt{\vglue -13truept{\logomed A\kern -.37em{\logobig 
T}\kern -.38em G}\vss}\hss}
\break
{\small Volume \thevolumenumber\ (\thevolumeyear)
\startpage--\finishpage\nl
Published: \publishdate}

\vglue .25truein

{\parskip=0pt\leftskip 0pt plus
1fil\def\\{\par\smallskip}{\Large\bf\thetitle}\par\medskip} \vglue
0.05truein

{\parskip=0pt\leftskip 0pt plus 1fil\def\\{\par}{\sc\theauthors}
\par\medskip}%
 
\vglue 0.03truein 

{\small\leftskip 25truept\rightskip 25truept{\bf Abstract}\stdspace\theabstract

{\bf AMS Classification}\stdspace\theprimaryclass
\ifx\thesecondaryclass\relax\else; \thesecondaryclass\fi\par
{\bf Keywords}\stdspace \thekeywords\par}\vglue 7truept

}   

\ifplaintex
\font\phead=cmsl9 scaled 950
\font\pnum=cmbx10 scaled 913
\font\pfoot=cmsl9 scaled 950
\headline{\vbox to 0pt{\vskip -4.5mm\line{\small\phead\ifnum
\count0=\startpage ISSN 1472-2739 (on-line) 1472-2747 (printed)
\hfill {\pnum\folio}\else\ifodd\count0\def\\{ }%
\ifx\theshorttitle\relax\thetitle\else\theshorttitle\fi\hfill{\pnum\folio}
\else\def\\{ and }{\pnum\folio}\hfill\ifx\theshortauthors\relax\theauthors
\else\theshortauthors\fi\fi\fi}\vss}}
\footline{\vbox to 0pt{\vglue 0mm\line{\small\pfoot\ifnum\count0=\startpage
\copyright\ \gtp\hfill\else
\agt, Volume \thevolumenumber\ (\thevolumeyear)\hfill\fi}\vss}}
\else
\font=cmsl9 scaled 1050
\font\lnum=cmbx10 
\font=cmsl9 scaled 1050
\makeatletter
\def\@oddhead{{\small\lhead\ifnum\count0=\startpage ISSN 1472-2739 
(on-line) 1472-2747 (printed)\hfill {\lnum\number\count0}\else\ifodd\count0
\def\\{ }\ifx\theshorttitle\relax \thetitle \else\theshorttitle\fi\hfill
{\lnum\number\count0}\else\def\\{ and }{\lnum\number\count0}
\hfill\ifx\theshortauthors\relax 
\theauthors\else\theshortauthors\fi\fi\fi}}\def\@evenhead{\@oddhead}
\def\@oddfoot{\small\lfoot\ifnum\count0=\startpage\copyright\ \gtp\hfill\else
\agt, Volume \thevolumenumber\ (\thevolumeyear)\hfill\fi}
\def\@evenfoot{\@oddfoot}
\makeatother
\fi
\let\maketitlepage\makeagttitle
\let\makeshorttitle\maketitlepage
\let\maketitle\maketitlepage

\newwrite\gtoutfile
\long\gdef\makeheadfile{  
{\def\\{, }\def\s{ }
\immediate\openout\gtoutfile head.xxx
\immediate\write\gtoutfile{To: math@arxiv.org}
\immediate\write\gtoutfile{Subject: put OR rep NNNNN:ppppp}
\immediate\write\gtoutfile{--text follows this line--}
\immediate\write\gtoutfile{Proxy-for: \ifx\theasciiauthors\relax
\theauthors\else\theasciiauthors\fi\s<\ifx\theasciiemail\relax\theemail\else\theasciiemail\fi>}
\immediate\write\gtoutfile{\noexpand\\}
\immediate\write\gtoutfile{Authors: \ifx\theasciiauthors\relax
\theauthors\else\theasciiauthors\fi}
{\def\\{ }\immediate\write\gtoutfile{Title: \ifx\theasciititle\relax
\thetitle\else\theasciititle\fi}}
\immediate\write\gtoutfile{Subj-class: GT or SG, GR etc}
\immediate\write\gtoutfile{MSC-class: \theprimaryclass\ifx\thesecondaryclass\relax\else, \thesecondaryclass\fi}
\immediate\write\gtoutfile{Journal-ref: Algebr. Geom. Topol. \thevolumenumber\s
(\thevolumeyear) \startpage-\finishpage}
\immediate\write\gtoutfile{Comments: Published by Algebraic and
Geometric Topology at}
\immediate\write\gtoutfile{\s\s\s  http://www.maths.warwick.ac.uk/agt/AGTVol\thevolumenumber/agt-\thevolumenumber-\thepapernumber.abs.html}
\immediate\write\gtoutfile{\noexpand\\}
\immediate\write\gtoutfile{}
\ifx\theasciiabstract\relax
\immediate\write\gtoutfile{\theabstract}\else
\immediate\write\gtoutfile{\theasciiabstract}\fi
\immediate\write\gtoutfile{}
\immediate\write\gtoutfile{\noexpand\\}
\immediate\write\gtoutfile{}
\immediate\closeout\gtoutfile}}  

\def\maketitlepage{\makeagttitle\makeheadfile}
\let\makeshorttitle\maketitlepage
\let\maketitle\maketitlepage

\lognumber{9}
\volumenumber{2}
\volumeyear{2002}
\papernumber{9}
\published{26 March 2002}
\pagenumbers{171}{217}
\received{3 December 2001}
\accepted{8 March 2002}

\usepackage{amsmath,amssymb}

\newtheorem{lem}{Lemma}[section]
\newtheorem{defi}[lem]{Definition}
\newtheorem{prop}[lem]{Proposition}
\newtheorem{rema}[lem]{Remark}
\newtheorem{theo}[lem]{Theorem}
\newtheorem{coro}[lem]{Corollary}
\newtheorem{exam}[lem]{Example}

\begin{document}
\title{Controlled connectivity of closed 1-forms}
\author{D. Sch\"utz}
\asciiauthors{D Schuetz}
\coverauthors{D Sch\noexpand\"utz}
\address{Department of Mathematics, University College Dublin\\ 
Belfield, Dublin 4, Ireland}
\email{dirk.schuetz@ucd.ie}

\begin{abstract}

We discuss controlled connectivity properties of closed 1-forms and
their cohomology classes and relate them to the simple homotopy type
of the Novikov complex. The degree of controlled connectivity of a
closed 1-form depends only on positive multiples of its cohomology
class and is related to the Bieri-Neumann-Strebel-Renz invariant. It
is also related to the Morse theory of closed 1-forms. Given a
controlled 0-connected cohomology class on a manifold $M$ with $n=\dim
M\geq 5$ we can realize it by a closed 1-form which is Morse without
critical points of index 0, 1, $n-1$ and $n$. If $n=\dim M\geq 6$ and
the cohomology class is controlled 1-connected we can approximately
realize any chain complex $D_\ast$ with the simple homotopy type of
the Novikov complex and with $D_i=0$ for $i\leq 1$ and $i\geq n-1$ as
the Novikov complex of a closed 1-form. This reduces the problem of
finding a closed 1-form with a minimal number of critical points to a
purely algebraic problem.
\end{abstract}

\asciiabstract{
We discuss controlled connectivity properties of closed 1-forms and
their cohomology classes and relate them to the simple homotopy type
of the Novikov complex. The degree of controlled connectivity of a
closed 1-form depends only on positive multiples of its cohomology
class and is related to the Bieri-Neumann-Strebel-Renz invariant. It
is also related to the Morse theory of closed 1-forms. Given a
controlled 0-connected cohomology class on a manifold M with n = dim M
> 4 we can realize it by a closed 1-form which is Morse without
critical points of index 0, 1, n-1 and n. If n = dim M > 5 and the
cohomology class is controlled 1-connected we can approximately
realize any chain complex D_* with the simple homotopy type of the
Novikov complex and with D_i=0 for i < 2 and i > n-2 as the Novikov
complex of a closed 1-form. This reduces the problem of finding a
closed 1-form with a minimal number of critical points to a purely
algebraic problem.}


\maketitle\def\\{\par}
\section{Introduction}
Given a finitely generated group $G$, Bieri, Neumann and Strebel
\cite{binest} and Bieri and Renz \cite{bieren} define subsets
$\Sigma^k(G)$ of equivalence classes of ${\rm
Hom}(G,\mathbb{R})-\{0\}$, where two homomorphisms $G\to \mathbb{R}$
are identified if they differ only by a positive multiple. These sets
reflect certain group theoretic properties of $G$ like finiteness
properties of kernels of homomorphisms to $\mathbb{R}$. In these
papers $\Sigma^k(G)$ is defined in terms of homological algebra but a
more topological approach is outlined as well. This topological
approach has become more important in recent years. Bieri and
Geoghegan \cite{biegeo} extend this theory to isometry actions of a
group $G$ on a ${\rm CAT}(0)$ space $M$. Although we will restrict
ourselves to the classical case, we will use this more modern approach
for our definitions. This way a property of a homomorphism $\chi:G\to \mathbb{R}$
being controlled $(k-1)$-connected ($CC^{k-1}$) is defined such that
$\chi$ being $CC^{k-1}$ is equivalent to $\pm[\chi]\in\Sigma^k(G)$. A
refinement which distinguishes between $\chi$ and $-\chi$ is also
discussed.\\
In the case where $G$ is the fundamental group of a closed connected smooth manifold $M$ the vector space
${\rm Hom}(G,\mathbb{R})$ can be identified with $H^1(M;\mathbb{R})$ via de Rham cohomology. Now the
controlled connectivity properties have applications in the Morse-Novikov theory of closed 1-forms.
Given a cohomology class $\alpha\in H^1(M;\mathbb{R})$ we can
represent it by a closed 1-form $\omega$ whose critical points are all
nondegenerate. We will call such 1-forms Morse forms. In particular there are
only finitely many critical points and every critical point has an index just as in ordinary Morse theory.
A natural question is whether there is a closed 1-form without critical points. This question was answered
by Latour in \cite{latour}. A similar problem is to find bounds for
the number of critical points of Morse forms representing $\alpha$ and
whether these bounds are exact. Special cases of this have been solved by Farber \cite{farbef} and
Pajitnov \cite{pajiol}.\\
To attack these problems one introduces the Novikov complex $C_\ast(\omega,
v)$ which first appeared in Novikov \cite{noviko}. This chain complex
is a free $\widehat{\mathbb{Z}G}_\chi$ complex generated by the
critical points of $\omega$ and graded by their indices. Here
$\widehat{\mathbb{Z}G}_\chi$ is a completion of the group ring $\mathbb{Z}G$ which depends on the homomorphism $\chi:G\to\mathbb{R}$ corresponding to the
cohomology class of $\omega$, see Section \ref{sec5} for details. 
To define the boundary in $C_\ast(\omega,v)$ one needs the
vector field $v$ to be gradient to $\omega$ and to satisfy
a transversality condition. This complex turns out to be simple chain
homotopy equivalent to
$\widehat{\mathbb{Z}G}_\chi\otimes_{\mathbb{Z}G}
C^\Delta_\ast(\tilde{M})$ where $\tilde{M}$ is the universal cover of $M$ and a
triangulation of $\tilde{M}$ is obtained by lifting a smooth triangulation of $M$. Therefore a closed 1-form
has to have at least as many critical points as any chain complex
$D_\ast$ has generators which is simple chain homotopic to
$\widehat{\mathbb{Z}G}_\chi\otimes_{\mathbb{Z}G}C^\Delta_\ast(\tilde{M})$.
Latour's theorem \cite[Th.$1'$]{latour} now reads as follows.
\begin{theo}\label{latint}
Let $M^n$ be a closed connected smooth manifold with $n\geq 6$ and $\alpha\in H^1(M;\mathbb{R})$.
Then $\alpha$ can be represented by a closed 1-form without critical points if and only if $\alpha$ is
$CC^1$, $\widehat{\mathbb{Z}G}_\chi\otimes_{\mathbb{Z}G}C^\Delta_\ast(\tilde{M})$ is acyclic and
$\tau(\widehat{\mathbb{Z}G}_\chi\otimes_{\mathbb{Z}G}C^\Delta_\ast(\tilde{M}))=0\in {\rm Wh}(G;\chi)$.
\end{theo}
Here ${\rm Wh}(G;\chi)$ is an appropriate quotient of
$K_1(\widehat{\mathbb{Z}G}_\chi)$. The condition that $\alpha$ be
$CC^1$ can be described as follows: a closed 1-form $\omega$
representing $\alpha$ pulls back to an exact form $df$ on the
universal cover.
For $\alpha$ to be $CC^1$ we require that for every
interval $(a,b)\subset\mathbb{R}$ there is a $\lambda\geq 0$ such that
every 0- or 1-sphere in $f^{-1}((a,b))$ bounds in 
$f^{-1}((a-\lambda,a+\lambda))$. Instead of $CC^1$, Latour \cite{latour} 
uses a stability condition on the homomorphism $\chi$ corresponding to
$\alpha$. We show in Section
\ref{sec4} that this is equivalent to our condition.\\
To prove this theorem one has to face the typical problems of the classical $h$- and $s$-cobordism
theorems. It turns out that the controlled connectivity conditions mentioned above are exactly what we need
for this. We get that $\alpha$ can be represented by a closed 1-form without critical points of index 0 and
$n=\dim M$ if and only if $\alpha$ is $CC^{-1}$. Of course this is equivalent to $\alpha$ being nonzero
and the corresponding fact that such a cohomology class can be represented without critical points of index
0 and $n$ has been known for a long time. If $n\geq 5$ removing critical points of index $0,1,n-1$ and $n$ is
equivalent to $CC^0$, see Section \ref{sec4}. Finally $CC^1$ allows us to perform the Whitney trick to
reduce the number of trajectories between critical points, provided $n\geq 6$. This is basically already
contained in Latour \cite[\S 4-5]{latour}, but we think that our approach is easier. Also the connection to
the Bieri-Neumann-Strebel-Renz theory in \cite{latour} is not mentioned.
Recently this connection was made more clear by Damian \cite{damian}, who also shows that the
condition $CC^1$ in Theorem \ref{latint} cannot be removed.\\
We deduce Latour's theorem by showing that for $\alpha$ $CC^1$ and $\dim M\geq 6$ we can realize a given
chain complex $D_\ast$ simple homotopy equivalent to $\widehat{\mathbb{Z}G}_\chi
\otimes_{\mathbb{Z}G}C^\Delta_\ast
(\tilde{M})$ approximately as the Novikov complex of a closed 1-form, provided $D_\ast$ is concentrated in
dimensions 2 to $n-2$. To be more precise, our main theorem is as follows.
\begin{theo}\label{introth}
Let $M^n$ be a closed connected smooth manifold with $n\geq 6$
and let $\alpha\in H^1(M;\mathbb{R})$ be
$CC^1$. Let $D_\ast$ be a finitely generated free based $\widehat{\mathbb{Z}G}_\chi$ complex with $D_i=0$ for
$i\leq 1$ and $i\geq n-1$ which is simple chain homotopy equivalent to $\widehat{\mathbb{Z}G}_\chi
\otimes_{\mathbb{Z}G}C^\Delta_\ast(\tilde{M})$. Given $L<0$ there is a Morse form $\omega$ representing
$\alpha$, a transverse $\omega$-gradient $v$ and a simple chain isomorphism $\varphi:D_\ast\to
C_\ast(\omega,v)$ where each $\varphi_i$ is of the form $I-A_i$ with $\|A_i\|<\exp L$.
\end{theo}
The negative real number $L$ comes from the fact that we do not actually realize the complex $D_\ast$
perfectly, but we can only approximate it arbitrarily closely.\\
A similar theorem has been proven by Pajitnov \cite[Th.0.12]{pajisu} in the case of a circle valued
Morse function $f:M\to S^1$. The condition $CC^1$ is replaced there by the condition that $\ker (f_\#:\pi_1
(M)\to\mathbb{Z})$ is finitely presented. This is in fact equivalent to $CC^1$ for rational closed 1-forms,
i.e.\ pullbacks of circle valued functions. See Theorem \ref{pajith} for a comparison to Pajitnov's
theorem.\\
In the exact case a similar theorem has been shown by Sharko \cite{sharko} which is in the same way a
generalization of the $s$-cobordism theorem as Theorem \ref{introth} is a generalization of Latour's
theorem.\\
Using Theorem \ref{introth} it is now easy to see that under the conditions that $\alpha$ is $CC^1$ and
$n\geq 6$ the minimal number of critical points of a closed 1-form within the cohomology class $\alpha$
is equal to the minimal number of generators of a chain complex $D_\ast$ of the simple homotopy type of
the Novikov complex. Thus the problem is reduced to a purely algebraic problem involving the Novikov ring
$\widehat{\mathbb{Z}G}_\chi$. Using the work of Farber and Ranicki \cite{farran} and Farber \cite{farber}
this problem can also be shifted to a different ring, a certain noncommutative localization of the group
ring, see Theorem \ref{localv} for more details.\\
As an application of Theorem \ref{introth} we can approximately
predescribe the torsion of a natural chain homotopy equivalence
$\varphi_v:\widehat{\mathbb{Z}G}_\chi\otimes C_\ast^\Delta(\tilde{M})
\to C_\ast(\omega,v)$ in $K_1(\widehat{\mathbb{Z}G}_\chi)/\langle [\pm
g]\,|\,g\in G\rangle$. The result we obtain is the following
\begin{theo}\label{introt2}
Let $G$ be a finitely presented group, $\chi:G\to\mathbb{R}$ be
$CC^1$, $b\in\widehat{\mathbb{Z}G}_\chi$ satisfy $\|b\|<1$ and
$\varepsilon>0$. Then for any closed connected smooth manifold $M$
with $\pi_1(M)=G$ and $\dim M\geq 6$ there is a Morse form $\omega$
realizing $\chi$, a transverse $\omega$-gradient $v$ and a
$b'\in\widehat{\mathbb{Z}G}_\chi$ with $\|b-b'\|<\varepsilon$ such
that $\tau(\varphi_v)=\tau(1-b')\in
K_1(\widehat{\mathbb{Z}G}_\chi)/\langle[\pm g]\,|\,g\in G\rangle$.
\end{theo}
By \cite[Th.1.1]{schue2} $\tau(\varphi_v)$
detects the zeta function of $-v$, a geometrically defined object
carrying information about the closed orbit structure of
$-v$. Therefore Theorem \ref{introt2} allows us to realize vector
fields whose zeta function is arbitrarily close to a predescribed
possible zeta function.\\
To prove Theorem \ref{introth} we have to realize certain elementary steps between simple chain homotopic
complexes for the geometric Novikov complexes.
The techniques of cancelling critical points, adding critical points and approximating an elementary change
of basis are all contained in Milnor \cite{milnhc}, but we have to make minor adjustments to be able to
use these methods in our situation. Most of these techniques in Milnor \cite{milnhc} are technically
quite involved, in order to not get hung up in technical difficulties
we mainly just write down the changes that
need to be done in the original proofs of \cite{milnhc}.\\
The results above suggest that vanishing of Novikov homology groups is related to controlled connectivity
conditions in general. To make this more precise one has to introduce a weaker notion called controlled
acyclicity. The precise relation can be found in Bieri \cite{bieri} or
Bieri and Geoghegan \cite{biege2}, but we discuss these results in
Section \ref{sec8} for the sake of completeness.\\
First results on the chain homotopy type of the Novikov complex were 
already announced in Novikov \cite{noviko}, but detailed proofs did not 
appear until much later, see Latour \cite{latour} or Pajitnov
\cite{pajito}. Easier proofs have since then appeared which are based on concrete chain homotopy
equivalences, but they are scattered through the literature and are not very well connected to each other.
In Appendix \ref{apchain} we describe some of these equivalences and show how they are related to each
other.\\
I would like to thank Ross Geoghegan for suggesting this topic and for
several valuable discussions. I would also like to thank Andrew
Ranicki for inviting me to Edinburgh where parts of this paper were
written. The author was supported by the EU under the TMR network FMRX
CT-97-0107 ``Algebraic K-Theory, Linear Algebraic Groups and Related
Structures''.
\subsection*{Notation}
Given a closed 1-form $\omega$ on a closed connected smooth manifold $M$ we denote the cohomology class
by $[\omega]\in H^1(M;\mathbb{R})$. A cohomology class $\alpha\in H^1(M;\mathbb{R})$ induces a homomorphism
$\chi=\chi_\alpha:\pi_1(M)\to\mathbb{R}$. We set $G=\pi_1(M)$. For a given closed 1-form $\omega$ there is a
minimal covering space such that $\omega$ pulls back to an exact form, namely the one corresponding to
$\ker \chi_{[\omega]}$. We denote it by $\rho:M_{[\omega]}\to M$. The universal covering space is denoted by
$\rho:\tilde{M}\to M$.
Given a vector field $v$ on $M$, we can lift it to covering spaces of $M$. We denote the lifting to
$\tilde{M}$ by $\tilde{v}$ and the lifting to $M_{[\omega]}$ by $\bar{v}$. If the critical points of $\omega$
are nondegenerate, we say $\omega$ is a Morse form. The set of critical points is denoted by
${\rm crit}\,\omega$.\\
Given a smooth function $f:N\to\mathbb{R}$ on a smooth manifold $N$ with nondegenerate critical points only
we define an $f$-gradient as in Milnor \cite[Df.3.1]{milnhc}, i.e.\ we have
\begin{enumerate}
\item $df(v)>0$ outside of critical points
\item if $p$ is a critical point of $f$, there is a neighborhood of $p$ such that $f=f(p)-\sum_{j=1}^i x_j^2
+\sum_{j=i+1}^n x_j^2$. In these coordinates we require $v=(-x_1,\ldots,-x_i,x_{i+1},\ldots,x_n)$.
\end{enumerate}
This notion of gradient extends in the obvious way to Morse forms. It is more restrictive than e.g.\ Pajitnov
\cite{pajiov} or \cite{schuet}, but is used to avoid further technicalities in cancelling critical
points.\\
Choose a Riemannian metric on $N$. If $p$ is a critical point and $\delta>0$, let $B_\delta(p)$, resp.\
$D_\delta(p)$, be the image of the Euclidean open, resp.\ closed, ball of radius $\delta$ under the
exponential map. Here $\delta$ is understood to be so small that $\exp$ restricts to a diffeomorphism of these
balls and so that for different critical points $p,\,q$ we get $D_\delta(p)\cap D_\delta(q)=
\emptyset$.\\
If $\Phi$ denotes the flow of an $f$-gradient $v$, we set:
\begin{eqnarray*}
W^s(p,v)&=&\{x\in N\,|\,\lim_{t\to\infty}\Phi(x,t)=p\}\\
W^u(p,v)&=&\{x\in N\,|\,\lim_{t\to-\infty}\Phi(x,t)=p\}\\
B_\delta(p,v)&=&\{x\in N\,|\,\exists t\geq 0 \hspace{0.4cm} \Phi(x,t)\in B_\delta(p)\}\\
D_\delta(p,v)&=&\{x\in N\,|\,\exists t\geq 0 \hspace{0.4cm} \Phi(x,t)\in D_\delta(p)\}
\end{eqnarray*}
The set $W^s(p,v)$ is called the stable and $W^u(p,v)$ the unstable manifold at $p$.
Notice that $W^s(p,v)\subset B_\delta(p,v)\subset D_\delta(p,v)$, $W^u(p,v)\subset B_\delta(p,-v)
\subset D_\delta(p,-v)$ and the $B_\delta$ sets are open. The sets $D_\delta(p,v)$ do not have to be closed
as other critical points might be in their closure.\\
A gradient $v$ is called transverse, if all stable and unstable manifolds intersect transversely. The set of
transverse gradients is generic, see Pajitnov \cite[\S 5]{pajiov}.\\
Let $R$ be a ring with unit and $\eta:\mathbb{Z}G\to R$ a ring homomorphism. Then define
\[C_\ast^\Delta(M;R)=R\otimes_{\mathbb{Z}G}C_\ast^\Delta(\tilde{M}) \mbox{ and } C_\Delta^\ast(M;R)=
{\rm Hom}_R(C_\ast^\Delta(M;R),R).\]
Notice that $C_\ast^\Delta(M;R)$ is a free left $R$ module and $C_\Delta^\ast(M;R)$ a free right
$R$ module. Furthermore we denote the homology and cohomology by $H_\ast(M;R)$ and $H^\ast(M;R)$.
\section{Controlled connectivity}
Let $k$ be a nonnegative integer and $G$ a group of type $F_k$, i.e.\ there exists a $K(G,1)$ CW-complex
with finite $k$-skeleton. Given a homomorphism $\chi:G\to\mathbb{R}$ we want to define statements "$\chi$
is controlled $(k-1)$-connected" and "$\chi$ is controlled $(k-1)$-connected over $\pm\infty$". To do this
let $X$ be the $k$-skeleton of the universal cover of a $K(G,1)$ CW-complex with finite $k$-skeleton. Then
$X$ is $(k-1)$-connected and $G$ acts freely and cocompactly on $X$ by covering translations. The
homomorphism $\chi$ induces an action of $G$ on $\mathbb{R}$ by translations, i.e.\ for $r\in\mathbb{R}$
we set $g\cdot r=r+\chi(g)$. An equivariant function $h:X\to\mathbb{R}$ is called a \em control function for
$\chi$. \em They exist because $G$ acts freely on $X$ and $\mathbb{R}$ is contractible. For $s\in\mathbb{R}$
and $r\geq 0$ denote $X_{s,r}(h)=\{x\in X\,|\,s-r\leq h(x)\leq s+r\}$. We will write $X_{s,r}$ if the
control function is clear.
\begin{defi}\em
The homomorphism $\chi:G\to\mathbb{R}$ is called \em controlled $(k-1)$-connected $(CC^{k-1})$, \em if for
every $r>0$ and $p\leq k-1$ there is a $\lambda\geq 0$ such that for every $s\in\mathbb{R}$ every $g:S^p
\to X_{s,r}$ extends to $\bar{g}:D^{p+1}\to X_{s,r+\lambda}$.
\end{defi}
This definition uses a choice of $X$ and $h$, but it turns out that controlled connectivity is a property
of $G$ and $\chi$ alone. To see that it does not depend on $h$ we have the following
\begin{lem}
Let $h_1,h_2:X\to\mathbb{R}$ be two control functions of $\chi$. Then there is a $t\geq0$ such that
$X_{s,r}(h_1)\subset X_{s,r+t}(h_2)\subset X_{s,r+2t}(h_1)$ for every $s\in\mathbb{R}$, $r\geq0$.
\end{lem}
\begin{proof}
Choose $t=\sup\{|h_1(x)-h_2(x)|\,|\,x\in X\}$, which is finite by cocompactness.
\end{proof}
\begin{lem}
The condition $CC^{k-1}$ does not depend on $X$.
\end{lem}
\begin{proof}
Let $Y_1$, $Y_2$ be two $K(G,1)$ CW-complexes with finite $k$-skeleton. Let $\alpha:Y_1\to Y_2$ and
$\beta:Y_2\to Y_1$ be cellular homotopy equivalences mutually inverse to each other. For $i=1,2$ let $X_i$
be the $k$-skeleton of the universal cover of $Y_i$. $\alpha$ and $\beta$ lift to maps $\tilde{\alpha}:
X_1\to X_2$ and $\tilde{\beta}:X_2\to X_1$ and we get a homotopy between $\tilde{\alpha}\circ\tilde{\beta}
|_{X_2^{(k-1)}}$ and the inclusion $X^{(k-1)}_2\subset X_2$, where $X^{(k-1)}_2$ denotes the $(k-1)$-skeleton.
Given a control function $h:X_2\to\mathbb{R}$ we get that $h\circ\tilde{\alpha}:X_1\to\mathbb{R}$ is also
a control function. Now $\tilde{\alpha}$ induces a map $(X_1)_{s,r}(h\circ\tilde{\alpha})\to(X_2)_{s,r}(h)$.
There is also a $t\geq0$ such that $\tilde{\beta}$ induces a map $(X_2)_{s,r}(h)\to (X_1)_{s,r+t}(h\circ
\tilde{\alpha})$ and we get a diagram
\[
\begin{array}{ccccc}
(X_1)_{s,r}(h\circ\tilde{\alpha})&\longrightarrow&(X_2)_{s,r}(h)&\longrightarrow&(X_1)_{s,r+t}(h\circ
\tilde{\alpha})\\[0.3cm]
\big\downarrow& &\big\downarrow& &\big\downarrow\\[0.3cm]
(X_1)_{s,r+\lambda}(h\circ\tilde{\alpha})&\longrightarrow&(X_2)_{s,r+\lambda}(h)&\longrightarrow&
(X_1)_{s,r+t+\lambda}(h\circ\tilde{\alpha})
\end{array}
\]
It follows that $\chi$ being $CC^{k-1}$ with respect to $X_2$ implies $\chi$ being $CC^{k-1}$ with
respect to $X_1$.
\end{proof}
It is clear that we can attach cells of dimension $\geq k+1$ to $X$ and still use $X$ to check for
$CC^{k-1}$. We also have $\chi$ is $CC^{k-1}$ if and only if $r\cdot\chi$ is $CC^{k-1}$ for
$r\not=0$.\\
Let us look at the case $k=0$. A $(-1)$-connected space is a nonempty space. Given a homomorphism
$\chi:G\to\mathbb{R}$, let us check for $CC^{-1}$. Choose $X$ and $h$. For $r>0$ we need $\lambda\geq 0$
such that for every $s\in\mathbb{R}$ the empty map $\emptyset\to X_{s,r}$ extends to $\bar{g}:\{\ast\}\to
X_{s,r+\lambda}$. So we need a $\lambda$ such that $X_{s,r+\lambda}$ is nonempty for all $s\in\mathbb{R}$.
This is clearly equivalent to $\chi$ being a nonzero homomorphism. A
nonzero homomorphism $\chi: G\to\mathbb{R}$ is also called
a \em character. \em\\
In the case where ${\rm im}\,\chi$ is infinite cyclic, $CC^0$ is equivalent to $\ker \chi$ being finitely
generated and $CC^1$ is equivalent to $\ker \chi$ being finitely presented. This follows from Brown
\cite[Th.2.2,Th.3.2]{brown} or Bieri and Geoghegan \cite[Th.A]{biegeo}.
\subsection*{Controlled connectivity over end points}
To draw a closer connection to the work of Bieri, Neumann and Strebel \cite{binest} and Bieri and Renz
\cite{bieren} let us define controlled connectivity over end points of $\mathbb{R}$. Let $X$ and $h$ be as
before. For $s\in\mathbb{R}$ define $X_s=\{ x\in X\,|\,h(x)\leq s\}$.
\begin{defi}\em
Let $\chi:G\to\mathbb{R}$ be a homomorphism. Then
\begin{enumerate}
\item $\chi$ is called \em controlled $(k-1)$-connected ($CC^{k-1}$) over
$-\infty$, \em if for every $s\in\mathbb{R}$ and $p\leq k-1$ there is a $\lambda(s)\geq 0$ such that every
map $g:S^p\to X_s$ extends to a map $\bar{g}:D^{p+1}\to X_{s+\lambda(s)}$ and $s+\lambda(s)\to -\infty$ as
$s\to -\infty$.
\item $\chi$ is called  \em controlled $(k-1)$-connected ($CC^{k-1}$) over
$+\infty$, \em if $-\chi$ is $CC^{k-1}$ over $-\infty$.
\end{enumerate}
\end{defi}
As before we get that these conditions only depend on $G$ and $\chi$, in fact they only depend on positive
multiples of $\chi$.\\
It is shown in Bieri and Geoghegan \cite[Th.H]{biegeo} that $\chi:G\to\mathbb{R}$ being $CC^{k-1}$ is
equivalent to $\chi$ being $CC^{k-1}$ at $-\infty$ and $+\infty$. This is also contained in Renz \cite{renz}.
Furthermore $\chi$ being $CC^{k-1}$
at $-\infty$ corresponds to $[\chi]\in \Sigma^k(G)$, the homotopical geometric invariant of Bieri and
Renz \cite[\S 6]{bieren}.
\subsection*{Cohomology classes and manifolds}
Now let $M$ be a closed connected smooth manifold and let $G=\pi_1(M)$. By de Rham's theorem we have
${\rm Hom}(G,\mathbb{R})=H^1(M;\mathbb{R})$ and we can represent cohomology classes by closed 1-forms
$\omega$. Now $\omega$ pulls back to an exact form on $\tilde{M}$, i.e.\ $\rho^\ast\omega=df$ with $f:
\tilde{M}\to \mathbb{R}$ smooth.
\begin{lem}\label{equimap}
The map $f:\tilde{M}\to\mathbb{R}$ is equivariant.
\end{lem}
\begin{proof}
Let $x\in\tilde{M}$, $g\in G$ and $\tilde{\gamma}$ a path from $x$ to $gx$. Then $\rho\circ\tilde{\gamma}$
represents the conjugacy class of $g\in G$ and we have
\[\chi(g)=\int\limits_{\rho_\ast\tilde{\gamma}}\omega=\int\limits_{\tilde{\gamma}}\rho^\ast\omega=
\int\limits_{\tilde{\gamma}}df=f(gx)-f(x),\]
so $f(gx)=f(x)+\chi(g)=g\cdot f(x)$.
\end{proof}
Therefore we can check for the controlled connectivity of $\chi$ by looking at a closed 1-form $\omega$
which represents $\chi$ and use the pullback $f$ as control function. Of course we need $\tilde{M}$ to
be $(k-1)$-connected to ask for $CC^{k-1}$, but we
can always check for controlled connectivity up to $CC^1$. In the special case of an aspherical $M$ on the
other hand we can check for $CC^{k-1}$ for any $k$. We will say $\alpha\in H^1(M;\mathbb{R})$ is $CC^{k-1}$,
if the corresponding homomorphism is. A control function of $\alpha$ will always refer to the pullback
of a closed 1-form representing $\alpha$.\\
Now assume that $\chi$ can be represented by a nonsingular closed 1-form $\omega$. Then $f:\tilde{M}\to
\mathbb{R}$ is a submersion. An $\omega$-gradient $v$ lifts to an $f$-gradient $\tilde{v}$ and we can use the
flowlines of $\tilde{v}$ to push every map $\bar{g}:D^p\to\tilde{M}$ into the subspace $X_{s,r}$.
So we can arrange $CC^{k-1}$ as long as $\tilde{M}$ is $(k-1)$-connected.\\
If we represent $\alpha$ by an arbitrary Morse form $\omega$ the critical points will represent an obstacle
to this approach. But if there exist no critical points of index less than $k$, a generic map $\bar{g}:
D^p\to\tilde{M}$ with $p<k$ will miss the unstable manifolds of the critical points of $f$ and we can use
the negative flow to get a map $\bar{g}_r:D^p\to X_s$ homotopic to $\bar{g}$ for every
$r\in\mathbb{R}$. So given a Morse form with no critical points of
index $<k$ and $>n-k$ we again get $CC^{k-1}$ as long as $\tilde{M}$
is $(k-1)$-connected.
\section{Changing a closed 1-form within a cohomology class}
The purpose of this section is to provide tools to modify a Morse form within its cohomology class. We need
to move the critical values of the control function in a useful way. This is achieved by starting with a
Morse form $\omega$ and a transverse $\omega$-gradient $v$ and modifying $\omega$ to a cohomologous form
$\omega'$ which agrees with $\omega$ near the critical points and such that $v$ is also an
$\omega'$-gradient. Then we need a tool to cancel critical points of $\omega$ in a nice geometric situation.
Both tools are described in Milnor \cite{milnhc}, but we need to sharpen the results to apply them to
irrational Morse forms, i.e.\ where the action induced by the form is not discrete. Compare also Latour
\cite[\S 3]{latour}.
\begin{lem}\label{movfun}
Let $N$ be a smooth manifold, $f:N\to\mathbb{R}$ a smooth function with nondegenerate critical points only
and $v$ an $f$-gradient. Let $p$ be a critical point of $f$, $\delta>0$ and $a<b$ such that $f(B_\delta(p))
\subset (a,b)$ and
$(D_\delta(p,v)\cup D_\delta(p,-v))\cap f^{-1}([a,b])$ contains no critical points except $p$. Then given
$c\in (a,b)$ there is a $g:N\to\mathbb{R}$ which agrees with $f$ outside of
$(D_\delta(p,v)\cup D_\delta(p,-v))\cap f^{-1}([a,b])$ such that $g(p)=c$ and $v$ is a $g$-gradient.
\end{lem}
\begin{proof}
Let $W=(D_\delta(p,v)\cup D_\delta(p,-v))\cap f^{-1}([a,b])$, $V=f^{-1}(\{a\})\cap W$ and $0<\delta_1<
\delta_2<\delta$. Define $\mu:V\to [0,1]$ to be 0 on $D_{\delta_1}(p,v)\cap V$ and bigger than $\frac{1}{2}$
on $V-D_{\delta_2}(p,v)$. Extend $\mu$ to $\bar{\mu}:W\to [0,1]$ by setting it constant on trajectories.
Now define $G:[a,b]\times[0,1]\to [a,b]$ with the properties
\begin{enumerate}
\item $\frac{\partial G}{\partial x}(x,y)>0$ and $G(x,y)$ increases from $a$ to $b$ as $x$ increases from
$a$ to $b$.
\item $G(f(p),0)=c$ and $\frac{\partial G}{\partial x}(x,0)=1$ for $x$ in a neighborhood of $f(p)$.
\item $G(x,y)=x$ for all $x$ if $y>\frac{1}{2}$ and for $x$ near $0$ and $1$ for all $y$.
\end{enumerate}
Now $g:W\to[a,b]$ defined by $g(q)=G(f(q),\bar{\mu}(q))$ extends to the desired function as in Milnor
\cite[Th.4.1]{milnhc}.
\end{proof}
Let $\omega$ be a Morse form, $v$ a transverse $\omega$-gradient and $f:M_\omega\to\mathbb{R}$ satisfy
$df=\rho^\ast\omega$. If $p\in M$ is a critical point of $\omega$, it lifts to a critical point $\bar{p}
\in M_\omega$ of $f$. Let $a<f(p)<b$ such that $A=(W^s(\bar{p})\cup W^u(\bar{p}))\cap f^{-1}([a,b])$ is a
positive distance away from all other critical points. Then $A$ is a compact set and since $v$ is transverse
we get that $A$ is disjoint from all translations of $A$ in $M_\omega$. Then there is a $\delta>0$ such that
this is also true for $(D_\delta(\bar{p},\bar{v})\cup D_\delta(\bar{p},-\bar{v}))\cap f^{-1}([a,b])$. So
we can apply Lemma \ref{movfun} equivariantly on $M_\omega$ to get the following.
\begin{lem}\label{movfor}
Let $\omega$ be a Morse form, $v$ a transverse $\omega$-gradient and $f:M_\omega\to\mathbb{R}$ the pullback
of $\omega$. If $\bar{p}\in M_\omega$ is a critical point of $f$ and $a<f(\bar{p})<b$ such that the closure
of $A=(W^s(\bar{p})\cup W^u(\bar{p}))\cap f^{-1}([a,b])$ contains no other critical points, then given
$c\in (a,b)$ and a neighborhood $U$ of $A$, there exists a Morse form $\omega'$ cohomologous to $\omega$
such that $v$ is an $\omega'$-gradient and a pullback $f':M_\omega\to\mathbb{R}$ that agrees with $f$
outside the translates of $U$ and satisfies $f'(\bar{p})=c$.
\end{lem}
\subsection*{Cancellation of critical points}
Theorem 5.4 of Milnor \cite{milnhc} shows how to cancel two critical points of adjacent index if there are
no other critical points around and there is exactly one trajectory between them. To apply this to our
situation we have to modify the result so that the function will only be changed in a neighborhood of the
critical points and part of the stable manifolds. More precisely we have
\begin{lem}\label{cancfun}
Let $f:N\to\mathbb{R}$ be smooth with nondegenerate critical points only and $v$ a transverse $f$-gradient.
Let $p,\,q$ be critical points with ${\rm ind}\,p={\rm ind}\,q+1$. Assume there is exactly one trajectory $T$
of $-v$ from $p$ to $q$ and an $\varepsilon>0$ such that for any other trajectory of $-v$ starting at $p$ and
ending in a critical point $p'$, we have $f(p')<f(q)-\varepsilon$. Then there is an arbitrarily small
neighborhood $V$ of $(W^s(p)\cup\{q\})\cap f^{-1}([f(q),f(p)])$ and a smooth function $f':N\to\mathbb{R}$
which agrees with $f$ outside $V$ and has no critical points in $V$. Furthermore there is an $f'$-gradient
$v'$ which agrees with $v$ outside an arbitrarily small neighborhood of $T$.
\end{lem}
\begin{proof}
Using Lemma \ref{movfun} we can change $f$ near $D_\delta(p,v)$ such that there is no trajectory of $v$
starting at $q$ and ending at a critical point $q'\not= p$ with $f(q')>f(p)+\varepsilon$, i.e.\ we can get
the images of $p$ and $q$ arbitrarily close together. So it is good enough to look at neighborhoods of the
form $U_\delta=(B_\delta(q,-v)\cup B_\delta(p,v))\cap f^{-1}((f(q)-\varepsilon,f(p)+\varepsilon))$ with
$\delta>0$ satisfying $\delta<\varepsilon$.\\
We can assume the Preliminary Hypothesis 5.5 of Milnor \cite{milnhc}. Using the first assertions of the
proof of Milnor \cite[Th.5.4]{milnhc} we can alter the vector field $v$ in $U_{\frac{\delta}{2}}$ to a
vector field $v'$ such that every trajectory starting in $U_{\frac{\delta}{2}}\cap f^{-1}(\{f(q)-
\frac{\delta}{2}\})$ reaches $f^{-1}(\{f(p)+\frac{\delta}{2}\})$ and stays within $U_{\frac{\delta}{2}}$.
Since $v'$ agrees with $v$ outside $U_{\frac{\delta}{2}}$ we get that $U_\delta$ is invariant for
trajectories of $v'$ within $f^{-1}((f(q)-\varepsilon,f(p)+\varepsilon))$. The closure of
$U_{\frac{\delta}{2}}\cap f^{-1}(\{f(q)-\frac{\delta}{2}\})$ is compact, so we get a product neighborhood
$V\times [0,1]\subset U_\delta$ of $T$ where $V=V\times\{0\}=U_\delta\cap f^{-1}(\{f(q)-\frac{\delta}{2}\})$
and $V\times\{1\}\subset f^{-1}(\{f(p)+\frac{\delta}{2}\})$. After rescaling we can assume $f(q)-
\frac{\delta}{2}=0$ and $f(p)+\frac{\delta}{2}=1$.\\
Consider $V\times [0,1]$as a subset of $N$ and define $g:V\times [0,1]\to\mathbb{R}$ by
\[g(x,u)=\int\limits_0^u\lambda(x,t)\frac{\partial f}{\partial t}(x,t)+(1-\lambda(x,t))\frac{
\int_0^1\lambda(x,s)\frac{\partial f}{\partial t}(x,s)\,ds}{\int_0^11-\lambda(x,s)\,ds}\,dt\]
where $\lambda:V\times[0,1]\to[0,1]$ is a smooth function which is constant 1 outside of
$U_{\frac{3}{4}\delta}\cap V\times[0,1]$ and in a small neighborhood of $V\times \{0,1\}$ and 0 in
$U_{\frac{\delta}{2}}\cap V\times [\varepsilon',1-\varepsilon']$ for $\varepsilon'>0$ so small that the
function is 0 where $v'$ differs from $v$. Notice that $g$ is smooth even for $x\in V$ with $\lambda(x,s)=1$
for all $s$.\\
As in Milnor \cite[p.54]{milnhc} it follows that $g$ extends to $f':N\to\mathbb{R}$ with the required
properties.
\end{proof}
For a Morse form $\omega$ the covering space $M_\omega$ has $G/\ker \chi$ as covering transformation group
and so there is a well defined homomorphism $\bar{\chi}:G/\ker\chi\to\mathbb{R}$. The desired Lemma to cancel
critical points of a Morse form now reads as
\begin{lem}\label{cancfor}
Let $\omega$ be a Morse form on the closed manifold $M$, $v$ a transverse $\omega$-gradient, $p,\,q$
critical points with ${\rm ind}\,p={\rm ind}\,q+1$. Let $\bar{p},\,\bar{q}\in M_\omega$ be lifts of $p$ and
$q$ such that there is exactly one trajectory $T$ between $\bar{q}$ and $\bar{p}$ and that there are no
trajectories between translates $D\bar{q}$ and $\bar{p}$ with $\bar{\chi}(D)>0$. Then there is a Morse form
$\omega'$ cohomologous to $\omega$ such that ${\rm crit}\,\omega'={\rm crit}\,\omega-\{p,q\}$ and an
$\omega'$-gradient $v'$ which agrees with $v$ outside an arbitrarily small neighborhood of $\rho(T)$.
\end{lem}
\begin{proof}
Let $f:M_\omega\to\mathbb{R}$ satisfy $df=\rho^\ast\omega$. Use Lemma \ref{movfor} to move the images of
the lifts of all critical points other than $q$ of index less than ${\rm ind}\,p$ by $f(\bar{p})-f(\bar{q})$
into the negative direction. To do this start with critical points of index 0, then critical points of index
1 and so on. This way we obtain a Morse form $\omega''$ with $\rho^\ast\omega''=df''$ and the same set of
critical points as $\omega$ which still has $v$ as a gradient. But now there are no trajectories between
$\bar{p}$ and critical points in $(f'')^{-1}([f''(\bar{q}),f''(\bar{p})])$ other than $T$. By choosing the
neighborhood $U_\delta$ in Lemma \ref{cancfun} small enough, we get that all translates of $U_\delta$ in
$M_\omega$ are disjoint. Now use Lemma \ref{cancfun} equivariantly on $f''$.
\end{proof}
\section{Relations between cancellation of critical points and controlled connectivity}\label{sec4}
We show that controlled connectivity in low degrees of a cohomology class leads to the existence of a Morse
form without critical points of low indices. This way we recover some well known results of Latour
\cite[\S 4]{latour} in a slightly different setting.
\begin{prop}\label{no0n}
Let $\alpha\in H^1(M;\mathbb{R})$. Then the following are equivalent:
\begin{enumerate}
\item $\alpha\not=0$.
\item $\alpha$ is $CC^{-1}$.
\item $\alpha$ can be represented by a Morse form $\omega$ without critical points of index $0,n$.
\end{enumerate}
\end{prop}
\begin{proof}
(1) $\Leftrightarrow$ (2) is clear.\\
(1) $\Rightarrow$ (3) : Choose an arbitrary Morse form $\omega$ and a transverse $\omega$-gradient $v$ and let
$p$ be a critical point of index 0. Lift $p$ to a critical point $\bar{p}\in M_\alpha$ and choose an
$\varepsilon>0$ such that the component of $\bar{p}$ in $A=f^{-1}((-\infty,f(\bar{p}+\varepsilon])$ is just
a small disc. Since $\alpha\not=0$ there are other components in $A$ or otherwise $\bar{p}$ would be an
absolute minimum of $f:M_\alpha\to\mathbb{R}$.\\
Claim: There is a critical point $q$ of $\omega$ of index 1 and a lift $\bar{q}\in M_\alpha$ such that one
of the flowlines of $W^s(\bar{q},\bar{v})$ ends in $\bar{p}$ while the other does not.\\
Since we know that $M_\alpha$ is connected there is a path between $\bar{p}$ and a point of $A$ which lies in
a different component of $A$. This path sits inside of some set $A'=f^{-1}((-\infty,f(\bar{p})+R])$ for
some big enough $R$. But if there is no critical point $\bar{q}$ as in the claim in $f^{-1}([f(\bar{p}),
f(\bar{p})+R])$ the component of $\bar{p}$ in $A'$ remains isolated by ordinary Morse theory.\\
Now the other trajectory of $\bar{q}$ can flow to
\begin{enumerate}
\item $-\infty$
\item a critical point $\bar{p}'$ of index 0 with $f(\bar{p}')<f(\bar{p})$
\item a critical point $\bar{p}'$ of index 0 with $f(\bar{p}')>f(\bar{p})$
\item a critical point $\bar{p}'$ of index 0 with $f(\bar{p}')=f(\bar{p})$.
\end{enumerate}
In the cases (1) and (2) we can cancel $p$ with $q$ by Lemma \ref{cancfor}. In case (3) we cancel $q$ with
$\rho(\bar{p}')$. In case (4) note that $\bar{p}'$ cannot be a translate of $\bar{p}$ because we are in
$M_\alpha$, so we can push the image of $\bar{p}$ slightly to a bigger number by Lemma \ref{movfor} and
then cancel $p$ and $q$. A dual argument holds for critical points of index $n$.\\
(3) $\Rightarrow$ (2) : Let $\omega$ be a Morse form without critical points of index $0,n$, $v$ a transverse
$\omega$-gradient and $f:\tilde{M}\to\mathbb{R}$ a control function. We claim that given $r>0$ and $x\in
\mathbb{R}$ there is a map $\bar{g}:\{\ast\}\to f^{-1}((x-r,x+r))$, i.e.\ we choose $\lambda=0$.
We know that $\tilde{M}$ is nonempty so let $y\in\tilde{M}$. Since there are no critical points of index
0 and $n$, any neighborhood of $y$ contains a dense subset of points that do not lie in any stable or
unstable manifold by Sard's theorem. Choose such a point. Using the flow of $\tilde{v}$ we can flow this
point to a point $y'$ with $f(y')=x$ for any $x\in\mathbb{R}$.
\end{proof}
Let $\alpha\in H^1(M;\mathbb{R})$ and $f:\tilde{M}\to\mathbb{R}$ a control function of $\alpha$. If $t\in
\mathbb{R}$ is a regular value we define
\[\tilde{N}(f,t)=f^{-1}(\{t\}).\]
\begin{lem}\label{finitapr}
Let $\omega$ be a Morse form and $v$ a transverse $\omega$-gradient. Let $t$ be a regular value of
$f:\tilde{M}\to\mathbb{R}$ where $df=\rho^\ast\omega$. Let $t_0>0$ and $C$ a compact subset of $\tilde{N}
(f,t)$. Then $C$ intersects only finitely many unstable discs $W^u(\tilde{p},\tilde{v})$ with $\tilde{p}
\in f^{-1}([t-t_0,t])$.
\end{lem}
\begin{proof}
For $i\geq 0$ define $W^i=\bigcup W^u(\tilde{p},\tilde{v})\cap f^{-1}([t-t_0,t])$ where the union is taken
over all critical points $\tilde{p}\in f^{-1}([t-t_0,t])$ with ${\rm ind}\,\tilde{p}\geq i$. Then $W^i$ is
closed. To see this notice that we can change $f$ on $f^{-1}((-\infty,t])$ to a function $g$ such that
$W^i\subset g^{-1}([t-t_0,t])$ and $g$ has no critical points of index $\leq i-1$ in $g^{-1}([t-t_0,t])$
and $\tilde{v}$ is a $g$-gradient. This is done using Lemma \ref{movfun} on every critical point in
$f^{-1}([t-t_0,t])$ of index $\leq i-1$ unequivariantly. So if $x\in g^{-1}([t-t_0,t])-W^i$, then $x$ is on
a trajectory going all the way to $g^{-1}(\{t\})=\tilde{N}(f,t)$. By continuity points near $x$ do the
same.\\
Therefore $W^i\cap C$ is compact. Now assume that $C$ intersects infinitely many discs. Since $\omega$ has
only finitely many critical points, there is a critical point $q$ such that $C$ intersects infinitely many
translates of $W^u(\tilde{q},\tilde{v})$ in $W^i$, where $\tilde{q}$ is a lifting of $q$ with $f(\tilde{q})
\in(t-t_0,t)$. Choose a point $x_k$ for every such translate. Since $C$ is compact there is an
accumulation point $x\in C$. Choose a small neighborhood $U$ of $x$ that gets mapped homeomorphically into
$M$ under the covering projection. Then there are infinitely many points $y_k\in W^u(\tilde{q},\tilde{v})$
and pairwise different $g_k\in G$ such that $g_ky_k\in U$ and $\{\chi(g_k)\}$ is bounded. But the $y_k$ also
have to have an accumulation point since $W^u(\tilde{q},\tilde{v})\cap f^{-1}([f(\tilde{q}),f(\tilde{q})+
t_0])$ has compact closure by the well definedness of the Novikov complex. But this contradicts $g_ky_k\in
U$ for infinitely many $k$.
\end{proof}
\begin{prop}\label{cc0con}
Let $\alpha\in H^1(M;\mathbb{R})$. Assume that $\alpha\not=0$ and $\dim M\geq 3$. Then the following are
equivalent:
\begin{enumerate}
\item $\alpha$ is $CC^0$.
\item There is a control function $f$ of $\alpha$ without critical points of index $0,n$ and with connected
$\tilde{N}(f,t)\subset\tilde{M}$.
\item There is a control function $f$ of $\alpha$ with connected $\tilde{N}(f,t)\subset\tilde{M}$.
\end{enumerate}
\end{prop}
\begin{proof}
(1) $\Rightarrow$ (2) : Choose $\omega'$ without critical points of index $0,n$ by Proposition \ref{no0n}
and a transverse $\omega'$-gradient $v$. Let $\tilde{N}'=\tilde{N}(f',t)$ where $t$ is a regular value
of $\tilde{f}'$ with $df'=\rho^\ast\omega'$. Since $\alpha$ is $CC^0$ there is a $\lambda>0$ such that
any two points in $\tilde{N}'$ can be connected in $\tilde{W}:=(f')^{-1}((t-\lambda,t+\lambda))$. Use Lemma
\ref{movfor} to get a new Morse form $\omega$ and control function $f$ such that $f(\tilde{p})=f'(\tilde{p})
-\lambda$ for every critical point $\tilde{p}$ of index 1 and $f(\tilde{q})=f'(\tilde{q})+\lambda$ for every
critical point of index $n-1$. Notice that since there are no critical points of index 0 and $v$ is
transverse, every critical point of index 1 can be pushed arbitrarily far to the negative side, similar for
critical points of index $n-1$.\\
Let $\tilde{N}=\tilde{N}(f,t)$. We claim that $\tilde{N}$ is connected.\\
Let $x,y\in\tilde{N}$. Since there are no critical points of index $0$ and $n$ we can assume that $x$ and $y$
do not lie on any stable or unstable manifold of $\tilde{v}$. So there are points $x',y'\in\tilde{N}'$ and
paths from $x$ to $x'$ and $y$ to $y'$ using flowlines. But $x'$ and $y'$ can be connected in $\tilde{W}$.
Using transversality we can find a smooth path between $x'$ and $y'$ in $\tilde{W}$ that does not meet any
stable or unstable manifolds of critical points $\tilde{q}$ with $2\leq{\rm ind}\,\tilde{q}\leq n-2$,
the stable manifolds of critical points with index 1 and the unstable manifolds of critical points with
index $n-1$.\\
By following flowlines this path can be pulled back into $\tilde{N}$ giving a path in $\tilde{N}$ between
$x$ and $y$. Assume not: let $z$ be a point on the path that lies on a trajectory that does not intersect
$\tilde{N}$. Without loss of generality assume $f(z)<t$, so $z$ would have to flow into the positive
direction to reach $\tilde{N}$. That the trajectory does not intersect $\tilde{N}$ means it converges to a
critical point $\tilde{p}$ with $f(\tilde{p})<r$ and ${\rm ind}\,\tilde{p}=n-1$ by the transversality
properties of the path. Now $f'(z)<f'(\tilde{p})=f(\tilde{p})-\lambda<t-\lambda$, contradicting $z\in
\tilde{W}$. So every point on the path can flow into $\tilde{N}$ giving a path between $x$ and $y$.\\
(2) $\Rightarrow$ (3) is trivial.\\
(3) $\Rightarrow$ (1) : Let $\tilde{N}=\tilde{N}(f,t)$ be connected. Choose $\lambda\geq 0$ such that $f^{-1}
([0,\lambda])$ contains two copies of $\tilde{N}$, note that $g\tilde{N}$ is a copy of $\tilde{N}$ in
$\tilde{M}$ for every $g\in G$. Let $g:S^0\to f^{-1}((x-r,x+r))$ be a map. Since $\tilde{M}$ is connected, we
can extend $g$ to a map $g':D^1\to\tilde{M}$. If $g'(D^1)\subset f^{-1}((x-r-\lambda,x+r+\lambda))$, we are
done. If not observe that by the choice of $\lambda$ both $f^{-1}((x-r-\lambda,x-r])$ and $f^{-1}([x+r,x
+r+\lambda))$ contain a copy of $\tilde{N}$. Denote them by $\tilde{N}_-$ and $\tilde{N}_+$. We can arrange
that $g'(D^1)$ intersects $\tilde{N}_-\cup\tilde{N}_+$ transversely. Then $(g')^{-1}(\tilde{N}_-\cup
\tilde{N}_+)\subset D^1$ is a finite set. Order them as $-1<t_1<\ldots<t_j<1$. If $g'(t_i)\not=g'(t_{i+1})$,
then the restriction of $g'$ to $[t_i,t_{i+1}]$ is a path in $f^{-1}((x-r-\lambda,x+r+\lambda))$. If
$g'(t_i)=g'(t_{i+1})$ we can change $g'$ on $[t_i,t_{i+1}]$ to a path in $\tilde{N}_\mp$. This way we get an
extension $\bar{g}:D^1\to f^{-1}((x-r-\lambda,x+r+\lambda))$ of $g$.
\end{proof}
\begin{prop}\label{no1n-1}
Let $\alpha\in H^1(M;\mathbb{R})$. Assume that $\alpha\not= 0$ and $\dim M\geq 5$. Then $\alpha$ is $CC^0$
if and only if $\alpha$ can be represented by a Morse form $\omega$ without critical points of index $0,1,
n-1,n$.
\end{prop}
\begin{proof}
Assume $\alpha$ is $CC^0$. Choose a Morse form $\omega$ without critical points of index $0,n$ and such that
there is a regular value $t\in\mathbb{R}$ with $\tilde{N}(f,t)$ connected by Proposition \ref{cc0con}. Let $v$
be a transverse $\omega$-gradient.\\
Let $p\in M$ be a critical point of index 1 and choose a lift $\tilde{p}\in \tilde{M}$ with $f(\tilde{p})>t$.
Let $r>f(\tilde{p})$ such that $\tilde{N}(f,r)\cap W^u(\tilde{p},\tilde{v})$ is an $(n-2)$-sphere $S$.
Denote the piece of the unstable manifold with boundary $S$ by $B$. Choose
a small arc in $\tilde{N}(f,r)$ that intersects $S$ transversely in one point and so that the endpoints do
not lie in any unstable manifold. Both endpoints can then flow into the negative direction until they reach
$\tilde{N}(f,t)$. Since $\tilde{N}(f,t)$ is connected we can choose a path between them. Now we have a loop
in $f^{-1}([t,r])$ which intersects $S$ transversely in exactly one point. We want to flow this loop back
to $\tilde{N}(f,r)$. By transversality we can change the loop so it avoids stable manifolds of critical points
with index $\leq n-2$. But we can change $\omega$ by Lemma \ref{movfor} by increasing the value of critical
points of index $n-1$ by $(r-t)$. By abuse of notation denote the resulting Morse form still by $\omega$
and $f$ for the control function. Then the loop can flow back to $\tilde{N}(f,r)$. Since $\tilde{M}$ is
simply connected, the loop bounds in $\tilde{M}$. By transversality we can embed a disc $D^2$ that avoids
stable manifolds of critical points of index $\leq n-3$. We can also arrange that $D^2$ embeds into $M$, not
just in $\tilde{M}$. Notice that $\partial D^2\subset\tilde{N}(f,r)$ and intersects $S$ in exactly one
point. Choose $a\leq b$ such that $f(D^2)\subset [a,b]$. Use Lemma \ref{movfor} to increase the value of
critical points of index $n-1$ and $n-2$ by $(b-a)$. Denote the resulting Morse form again by $\omega$ and
the control function by $f$. Note that this can be done so that $\partial D^2$ and $S$ are still in
$\tilde{N}(f,r)$. We can assume that $b$ is a regular value. Now we can use the flow of $\tilde{v}$ to push
$D^2$ into $\tilde{N}(f,b)$. Denote the boundary of that disc by $S_1$. We have that $S_1$ intersects
$W^u(\tilde{p},\tilde{v})$ transversal in exactly one point, $S_1$ embeds into $M$ and $S_1$ bounds a disc
$D_1^2$ in $\tilde{N}(f,b)$. Since the vector field will be changed in a small neighborhood of $D^2_1$, we
need to make sure that $D^2_1$ is nice. Since $S_1$ is obtained from $\partial D^2$ by flowing we get a
2-dimensional surface $S^1\times I$ between $S_1$ and $\partial D^2$. Use transversality to modify $D^2_1$
such that it does not intersect any translates of that surface. We do not want to change the vector field
on $B$. Since $B$ is $(n-1)$-dimensional and $D_1^2\subset\tilde{N}(f,b)$, $D^2_1$ can intersect translates
of $B$ in finitely many circles. But whenever we have such a circle, we can change $D^2_1$ to remove the
intersection since the normal bundle of $B$ is trivial.\\
Now we can proceed as in Milnor \cite[p.105]{milnhc}. Insert two critical points
$\tilde{q}$, $\tilde{q}'$ of index 2 and 3 equivariantly near the right of $S_1$. Adjust $\tilde{v}$ to
$\tilde{v}'$ so that $W^s(\tilde{q},\tilde{v})\cap\tilde{N}(f,b)=S_1$. Then there is exactly one flowline
from $\tilde{q}$ to $\tilde{p}$ and all other trajectories from $\tilde{q}$ go to the left of $\tilde{p}$.
Hence we can cancel $\tilde{p}$ and $\tilde{q}$. This way we can trade all critical points of index 1 for
critical points of index 3. A dual argument works for critical points of index $n-1$.\\
Now assume we have a control function $f$ without critical points of index $0,1,n-1,n$. Given $g:S^0\to
f^{-1}((x-r,x+r))$ it extends to a map $D^1\to\tilde{M}$. By transversality we can change this map so that
it avoids stable and unstable manifolds in the interior of $D^1$. Then we can use the flow to push it into
$f^{-1}((x-r,x+r))$.
\end{proof}
\begin{prop}\label{cc1con}
Let $\alpha\in H^1(M;\mathbb{R})$. Assume that $\alpha\not=0$ and $\dim M\geq 5$. Then the following are
equivalent:
\begin{enumerate}
\item $\alpha$ is $CC^1$.
\item There is a control function $f$ of $\alpha$ without critical points of index $0,1,n-1,n$ and with
simply connected $\tilde{N}(f,t)\subset\tilde{M}$.
\item There is a control function $f$ of $\alpha$ with simply connected $\tilde{N}(f,t)\subset\tilde{M}$.
\end{enumerate}
\end{prop}
\begin{proof}
(1) $\Rightarrow$ (2) : The proof is analogous to the proof of Proposition \ref{cc0con}. Choose a Morse form
$\omega'$ representing $\alpha$ without critical points of index $0,1,n-1,n$ by Proposition \ref{no1n-1}, let
$v$ be a transverse $\omega'$-gradient and let $f':\tilde{M}\to\mathbb{R}$ be a control function.
Let $\tilde{N}'=\tilde{N}(f',t)$ where $t\in\mathbb{R}$ is a regular value. Since $\alpha$ is $CC^1$ there
is a $\lambda>0$ such that every loop in $\tilde{N}'$ bounds in $\tilde{W}:=(f')^{-1}((t-\lambda,
t+\lambda))$. Change $\omega'$ to a Morse form $\omega$ with control function $f$ such that $f(\tilde{p})=
f'(\tilde{p})-\lambda$ for critical points of index 2 and $f(\tilde{q})=f'(\tilde{q})+\lambda$ for critical
points of index $n-2$.\\
Let $\tilde{N}=\tilde{N}(f,t)$ and $\gamma$ a loop in $\tilde{N}$. Using transversality we can assume that
$\gamma$ does not intersect any stable or unstable manifolds of $\tilde{v}$. So we can use the flow of
$\tilde{v}$ to flow $\gamma$ into $\tilde{N}'$. This loop bounds in $\tilde{W}$. Choose the disc so that it
intersects stable and unstable manifolds transversely. This disc now flows back into $\tilde{N}$ as in the
proof of Proposition \ref{cc0con}. Therefore $\tilde{N}$ is simply connected.\\
(2) $\Rightarrow$ (3) is trivial.\\
(3) $\Rightarrow$ (1) : Let $\tilde{N}=\tilde{N}(f,t)$ be simply connected. Then $\alpha$ is $CC^0$ by Proposition
\ref{cc0con}. Choose $\lambda\geq 0$ such that $f^{-1}([0,\lambda])$ contains two copies of $\tilde{N}$.
Let $g:S^1\to f^{-1}((x-r,x+r))$ be a map. This extends to a map $g':D^2\to \tilde{M}$. Let $\tilde{N}_-$
be a copy of $\tilde{N}$ in $f^{-1}((x-r-\lambda,x-r])$ and $\tilde{N}_+$ a copy of $\tilde{N}$ in
$f^{-1}([x+r,x+r+\lambda))$. We can assume that $g'$ intersects $\tilde{N}_-\cup\tilde{N}_+$ transversely,
i.e.\ in a finite set of circles. Since these circles bound in $\tilde{N}_\mp$ we can change $g'$ away
from the boundary to a map $\bar{g}:D^2\to f^{-1}((x-r-\lambda,x+r+\lambda))$.
\end{proof}
\begin{rema}\em
Latour \cite{latour} uses a stability condition on $\pm\chi$ instead of 
$CC^1$ to obtain \cite[Prop.5.20]{latour} which is analogous to Proposition
\ref{cc1con}. It follows that $\chi$ being $CC^1$ and $\pm \chi$ being
stable in the sense of Latour \cite[\S 5]{latour} are equivalent. The condition
of only $\chi$ being stable corresponds to $\chi$ being $CC^1$ at $\infty$,
compare Section \ref{sec8} and \cite[Cor.5.16]{latour}.
\end{rema}
\section{The Novikov complex}\label{sec5}
Let $G$ be a group and $\chi:G\to\mathbb{R}$ be a homomorphism. We denote by $\widehat{\widehat{\mathbb{Z}G}}$
the abelian group of all functions $G\to\mathbb{Z}$. For $\lambda\in\widehat{\widehat{\mathbb{Z}G}}$
let supp $\lambda=\{g\in G\,|\,\lambda(g)\not=0\}$. Then we define
\[\widehat{\mathbb{Z}G}_\chi=\{\lambda\in\widehat{\widehat{\mathbb{Z}G}}\,|\,\forall r\in\mathbb{R}
\hspace{0.4cm}\#\,\mbox{supp }\lambda\cap\chi^{-1}([r,\infty))<\infty\}.\]
For $\lambda_1,\lambda_2\in\widehat{\mathbb{Z}G}_\chi$ we set $(\lambda_1\cdot\lambda_2)(g)=\sum
\limits_{h_1,h_2\in G\atop h_1h_2=g}\lambda_1(h_1)\lambda_2(h_2)$, then $\lambda_1\cdot\lambda_2$
is a well defined element of $\widehat{\mathbb{Z}G}_\chi$ and turns $\widehat{\mathbb{Z}G}_\chi$ into
a ring, the \em Novikov ring\em. It contains the usual group ring $\mathbb{Z}G$ as a subring and
we have $\mathbb{Z}G=\widehat{\mathbb{Z}G}_\chi$ if and only if $\chi$ is the zero homomorphism.
\begin{defi}\label{dnorm}\em The \em norm \em of $\lambda\in\widehat{\mathbb{Z}G}_\chi$ is defined to be
\[\|\lambda\|=\|\lambda\|_\chi=\inf\{t\in(0,\infty)|\mbox{ supp }\lambda
\subset\chi^{-1}((-\infty,\log t])\}.\]
\end{defi}
For $L\in\mathbb{R}$ define $p_L:\widehat{\mathbb{Z}G}_\chi\to\widehat{\mathbb{Z}G}_\chi$ by $p_L(\lambda)(g)
=\left\{\begin{array}{cl}\lambda(g)&\chi(g)\geq L\\0&\mbox{otherwise}\end{array}\right.$. Notice that $p_L$
factors through $\mathbb{Z}G$ and is a homomorphism of abelian groups, but not of rings. It also
extends to free $\widehat{\mathbb{Z}G}_\chi$ modules.\\
Given $a\in\widehat{\mathbb{Z}G}_\chi$ with $\|a\|<1$, the series $\sum_{k=0}^\infty a^k$ is a well defined
element of $\widehat{\mathbb{Z}G}_\chi$ and hence the inverse of $1-a$. Therefore $\{1-a\in\widehat{\mathbb{Z}
G}_\chi\,|\,\|a\|<1\}$ is a subgroup of the group of units. Let ${\rm Wh}(G;\chi)$ be the quotient of $K_1
(\widehat{\mathbb{Z}G}_\chi)$ by these units and units of the form $\pm g$ with $g\in G$.\\
Given a Morse form $\omega$ and a transverse $\omega$-gradient $v$ we can define the Novikov complex
$C_\ast(\omega,v)$ which is in each dimension $i$ a free $\widehat{\mathbb{Z}G}_\chi$ complex with one
generator for every critical point of index $i$. Here $\chi$ is the homomorphism induced by $\omega$. To
define the boundary homomorphism choose an orientation for the stable manifolds of every critical point.
Now coorient the unstable manifolds, i.e.\ choose an orientation of the normal bundle so that the
coorientation at $W^u(p,v)$ projects to the chosen orientation of $W^s(p,v)$ at $p$. If $p,\,q$ are critical
points with ${\rm ind}\,p={\rm ind}\,q+1=i$, then $W^s(p,v)\cap W^u(q,v)$ is 1-dimensional which means it
consists of isolated trajectories. Given a trajectory $T$ between $p$ and $q$ we want to define a sign
for $T$. If $x\in T$ let $X\in T_xM$ be a vector with $\omega(X)<0$. Also let $X_1,\ldots,X_{i-1}\in T_xM$
represent the coorientation of $W^u(q,v)$. If the projection of $X,X_1,\ldots,X_{i-1}$ into the
tangent space of $W^s(p,v)$ at $x$ represents the orientation of $W^s(p,v)$, set $\varepsilon(T)=1$, otherwise
set $\varepsilon(T)=-1$. Note that these projections do represent a basis for $T_xW^s(p,v)$ by the
transversality assumption.\\
Now lift the orientations to $\tilde{M}$ and choose for every critical point of $\omega$ exactly one lift
in $\tilde{M}$. For critical points $p,\,q$ with ${\rm ind}\,p={\rm ind}\,q+1$ define $[p:q]\in\widehat{
\mathbb{Z}G}_\chi$ by
\[[p:q]\,(g)=\sum\varepsilon(T)\]
where the sum is taken over the set of all trajectories between $\tilde{p}$ and $g\tilde{q}$, where
$\tilde{p}$ and $\tilde{q}$ are the chosen liftings of $p$ and $q$. Then define $\partial:C_\ast(\omega,v)
\to C_{\ast-1}(\omega,v)$ by
\[\partial(p)=\sum_{q,\,{\rm ind}\,q={\rm ind}\,p-1}[p:q]\,q.\]
That $[p:q]$ is indeed an element of $\widehat{\mathbb{Z}G}_\chi$ and $\partial^2=0$ is shown in the exact
case in Milnor \cite[\S 7]{milnhc}. The case of a circle valued Morse function can be reduced to the exact
case by inverse limit arguments, compare Pajitnov \cite{pajito} or Ranicki \cite{ranick}. Finally the
irrational case can be reduced to the rational case by approximation, see Pajitnov \cite{pajisp} or the
author \cite[\S 4.2]{schuet}.\\
Since $[p:q]$ depends on the gradient $v$, we also write $[p:q]_v$ when we deal with different
gradients.\\
The appendix describes simple chain homotopy equivalences $\varphi_v:C_\ast^\Delta(M;
\widehat{\mathbb{Z}G}_\chi)$ $\to C_\ast(\omega,v)$ and $\psi_{v_1,v_2}:C_\ast(\omega_1,v_1)\to
C_\ast(\omega_2,v_2)$ with $\psi_{v_1,v_2}\circ\varphi_{v_1}\simeq\varphi_{v_2}$, $\psi_{v_2,v_3}\circ
\psi_{v_1,v_2}\simeq\psi_{v_1,v_3}$ and $\psi_{v_1,v_1}\simeq{\rm id}$, where $\simeq$ means chain
homotopic and $\omega,\,\omega_1,\,\omega_2,$ $\omega_3$ are cohomologous.\\
To define the Novikov complex, we made a choice of liftings of critical points. Let $\mathcal{B}\subset
{\rm crit}(f)$ be this choice. Set $A_\mathcal{B}=\sup\{|f(\tilde{p})-f(\tilde{q})|\,|\,\tilde{p},\,
\tilde{q}\in\mathcal{B}\}$.
\begin{prop}\label{elmchan}
Let $\omega$ be a Morse form without critical points of index $0,1,n-1,n$, $v$ a transverse $\omega$-gradient,
$\mathcal{B}$ a choice of liftings of the critical points of $\omega$ and $p_1\not= p_2$ critical points of
$\omega$ having index $i$. Let $g\in G$ be such that $f(\tilde{p}_1)>f(\tilde{p}_2)+\chi(g)$, where
$\tilde{p}_1,\,\tilde{p}_2\in\mathcal{B}$ are liftings of $p_1,\,p_2$ and $L<\min\{0,\,\chi(g)\}$. Then
there is a transverse $\omega$-gradient $v'$ such that:
\begin{enumerate}
\item $p_L(\psi_{v',v}(p_1))=p_1+gp_2$.
\item $p_L(\psi_{v',v}(q))=q$ for critical points $q\not=p_1$.
\item $p_L([p:q]_{v'})=p_L([p:q]_v)$ for $p\not=p_1$, ${\rm ind}\,q={\rm ind}\,p-1$.
\item $p_L([p_1:q]_{v'})=p_L([p_1:q]_v)+gp_L([p_2:q]_v)$ for ${\rm ind}\,q=i-1$.
\end{enumerate}
\end{prop}
We can think of the statement as performing an elementary change of basis, but we can only approximate the
elementary change. The proof is based on Milnor \cite[Th.7.6]{milnhc}. The condition $f(\tilde{p}_1)>
f(\tilde{p}_2)+\chi(g)$ can always be achieved by changing $\omega$ using Lemma \ref{movfor}.
\begin{proof}
Choose a regular value $t_0$ with $f(\tilde{p_1})>t_0>f(\tilde{p_2})+\chi(g)$ and set $V_0=f^{-1}(\{t_0\})$.
We have $S_L:=W^s(\tilde{p}_1,\tilde{v})\cap V_0$ is $(i-1)$-dimensional and $S_R:=W^u(g\tilde{p}_2,
\tilde{v})\cap V_0$ is $(n-i-1)$-dimensional. Since there are no critical points of index $0,1,n-1,n$, both are
nonempty and $V_0$ is connected. Hence we can embed a path $\varphi:[0,3]\to V_0$ that intersects $S_L$
transversely at $\varphi(1)$, $S_R$ transversely at $\varphi(2)$ and that misses all other stable and
unstable manifolds. Using Milnor \cite[Lm.7.7]{milnhc} we get a nice product neighborhood $U$ of this arc.
By choosing it small enough and Lemma \ref{finitapr} we can assume that it misses the unstable and stable
manifolds of critical points
of $f$ other than $\tilde{p}_1$ and $g\tilde{p}_2$ in $f^{-1}((t_0-A_\mathcal{B}+L,t_0+A_\mathcal{B}-L))$.
To see this notice that for every critical point of $\omega$ there are only finitely many liftings in
$f^{-1}((t_0-A_\mathcal{B}+L,t_0+A_\mathcal{B}-L))$ whose stable or unstable manifolds can get close to
the arc. So they stay away a positive distance.\\
Now using the flow of $\tilde{v}$ we can find a small product neighborhood of $U$ in $\tilde{M}$ and change
the vector field $\tilde{v}$ to a vector field $\tilde{v}''$ equivariantly as in Milnor \cite[p.96]{milnhc}.
The stable and unstable manifolds of critical points of $f$ other than $\tilde{p}_1$ and $\tilde{p}_2$
do not get changed within a range of $\pm (A_\mathcal{B}-L)$. The $\omega$-gradient $v''$ need not be
transverse, but we can find a transverse $\omega$-gradient $v'$ as close as we like to $v$ in the
smooth topology. Choose one so close that the intersection numbers of stable and unstable manifolds
within the $\pm (A_\mathcal{B}-L)$ range are as with
$v''$. Then the properties (1)-(4) of $\psi_{v',v}$ follow by the definition of $\psi_{v',v}$ and the fact
that we can use $v''$ for $p_L([p:q]_{v'})$.
\end{proof}
For the next proposition the controlled 1-connectivity is crucial.
\begin{prop}\label{canceltraj}
Let $\omega$ be a Morse form without critical points of index $0,1,n-1,n$ which is $CC^1$ and $v$ a
transverse $\omega$-gradient. Assume that $n=\dim M\geq 6$. Let $q$ be a critical point of index $i$ with
$2\leq i\leq n-3$ and $p$ be a critical point of index $i+1$. Let $\tilde{p}$, $\tilde{q}$ be liftings of
$p$ and $q$ to $\tilde{M}$ such that there exist two trajectories $T_1,\,T_2$ between $\tilde{p}$ and
$\tilde{q}$ with $\varepsilon(T_1)=-\varepsilon(T_2)$ and there exist no trajectories between $\tilde{p}$
and $g\tilde{q}$ with $\chi(g)>0$. Let $L<0$. Then there is a Morse form $\omega'$
cohomologous to $\omega$ which agrees with $\omega$ at the common set of critical points, a transverse
$\omega'$-gradient $v'$ such that there are two less trajectories between $\tilde{p}$ and $\tilde{q}$,
no new trajectories between $\tilde{p}$ and $g\tilde{q}$ with $\chi(g)>L$ and we have $p_L\circ
\psi_{v',v}=p_L$ and $p_L([r:s]_{v'})=p_L([r:s]_v)$ for ${\rm ind}\,r={\rm ind}\,s+1$.
\end{prop}
\begin{proof}
Let us assume that $2\leq i\leq n-4$, if $i=n-3$, look at $-\omega$ and $-v$.\\
We can alter $\omega$ as in the proof of Proposition \ref{cc1con} such that there is a simply connected
$\tilde{N}(f,t)$. In the irrational case we can assume that $t$
satisfies $f(\tilde{q})>t>f(\tilde{p})$ and that $t$ is so close to $f(\tilde{q})$ such that $\tilde{N}
(f,t)\cap W^u(\tilde{q},\tilde{v})$ is a sphere of dimension $(n-i-1)$. In the rational case we can change
$\omega$ so that $f$ orders the critical points in $f^{-1}([t,t+\chi(g)])$, where $\chi(g)$ generates
${\rm im}\,\chi$
and then we also get a simply connected $\tilde{N}(f,t)$ with $f(\tilde{q})<t<f(\tilde{p})$ and $\tilde{N}
(f,t)\cap W^u(\tilde{q},\tilde{v})$ is a sphere of dimension $(n-i-1)$. We now want $\tilde{N}(f,t)\cap
W^s(\tilde{p},\tilde{v})$ to be a sphere of dimension $i$. Since there are no trajectories between
$\tilde{p}$ and $g\tilde{q}$ with $\chi(g)>0$ we can achieve this by changing $\omega$ to a Morse form
$\omega'$ such that $f'(\tilde{r})=f(\tilde{r})-(f(\tilde{p})-t)$ for every critical point $r\not= q$ with
${\rm ind}\,r\leq i$.\\
Let $\tilde{N}=\tilde{N}(f',t)$. Then $\tilde{N}\cap W^u(\tilde{q},\tilde{v})=S_R$ and $\tilde{N}\cap
W^s(\tilde{p},\tilde{v})=S_L$ are spheres.\\
We need that $\tilde{N}$ is still simply connected. But a loop in $\tilde{N}$ is homotopic to one that
can be flown into $\tilde{N}(f,t)$ since there are no critical points of index 0 and 1. Now this loop bounds
in the simply connected $\tilde{N}(f,t)$. But a generic 2-disc can flow back into $\tilde{N}$, since we only
moved critical points of index $\leq n-4$. This shows that $\tilde{N}$ is also simply connected.\\
We want to apply Milnor \cite[Th.6.6]{milnhc}. To see that $\tilde{N}-S_R$ is simply connected the same
argument as in Milnor \cite[p.72]{milnhc} works. Notice that the isotopy in \cite[Th.6.6]{milnhc} is fixed
outside a neighborhood of a 2-disc which bounds two arcs between $T_1\cap\tilde{N}$ and $T_2\cap\tilde{N}$.
By transversality arguments we can assume that this disc does not intersect any unstable manifolds
$W^u(\tilde{r},\tilde{v})$ and stable manifolds $W^s(\tilde{s},\tilde{v})$ for $\tilde{r},\,\tilde{s}\in
(f')^{-1}((t-A_\mathcal{B}+L,t+A_\mathcal{B}-L))$ with ${\rm ind}\,\tilde{r}\geq i+1$ and ${\rm ind}\,
\tilde{s}\leq i$. Here $A_\mathcal{B}$ is defined as before with respect to $f'$. We can also
assume the disc embeds into $M$. By choosing the neighborhood of the disc small enough we can change the
$\omega'$-gradient $v$ to an $\omega'$-gradient $v''$ with two fewer trajectories between $\tilde{p}$ and
$\tilde{q}$. Choose a transversal $\omega'$-gradient $v'$ so close to $v''$ such that the intersection
numbers in the $\pm(A_\mathcal{B}-L)$-range are the same. By the way the neighborhood of the disc was
chosen we now get $p_L\circ\partial'=p_L\circ\partial$ and $p_L\circ\psi_{v',v}=p_L$.
\end{proof}
\begin{prop}\label{insert}
Let $\omega$ be a Morse form without critical points of index $0,1,n-1,n$ and $v$ a transverse
$\omega$-gradient. Let $x\in M$ be a regular point, $i$ an integer with $2\leq i\leq n-3$ and $L<0$.
Given any neighborhood $U$ of $x$ there is a Morse form $\omega'$ and a transverse
$\omega'$-gradient $v'$ such that $\omega'$ agrees with $\omega$ outside $U$ and ${\rm crit}\,\omega'=
{\rm crit}\,\omega\cup\{p,q\}$ with $p,\,q\in U$ and ${\rm ind}\, p=i+1$, ${\rm ind}\,q=i$ such that
\begin{enumerate}
\item $[p:q]_{v'}=1-a$ with $\|a\|<1$.
\item $\max\{\|[p:q']_{v'}\|,\,\|[r:p]_{v'}\|,\,\|[p':q]_{v'}\|,\,\|[q:s]_{v'}\|\}<\exp L$, where
${\rm ind}\,q'=i$, ${\rm ind}\,p'=i+1$, ${\rm ind}\,r=i+2$, ${\rm ind}\,s=i-1$, $q'\not=q$ and $p'\not=p$.
\item $p_L([p':q']_{v'})=p_L([p':q']_{v})$ for $p'\not=p$ and $q'\not=
q$.
\item $p_L(\psi_{v,v'}(p'))=p'$ for $p'\in{\rm crit}\,\omega$.
\end{enumerate}
\end{prop}
\begin{proof}
Since there are no critical points of index $0$ and $n$ there is a $y\in U$ which does not lie on any stable
or unstable manifold. Let $\tilde{y}\in\tilde{M}$ be a lift of $y$. We can find a small neighborhood $V$ of
$y$ with $V\subset U$ such that $\tilde{V}\cap W^{s,u}(\tilde{r},\tilde{v})=\emptyset$, if $r$ is a
critical point of $\omega$ and $|f(\tilde{r})-f(\tilde{y})|<2A_\mathcal{B}-L$. Here $\tilde{V}$ is a lift
of $V$ with $\tilde{y}\in\tilde{V}$.\\
We can insert two critical points $p,\,q$ of adjacent indices as in Milnor \cite[p.105]{milnhc}. This way we
obtain a Morse form $\omega'$ and a $\omega'$-gradient $v''$. Choose a transverse $\omega'$-gradient
$v'$ so close to $v''$ the intersection numbers do not change in a range of $\pm(A_\mathcal{B}-L)$. Choose
the liftings for the basis of the corresponding Novikov complex in a translate of $\tilde{V}$ such that
$A_\mathcal{B}$ does
not increase by more than $|f(\tilde{p})-f(\tilde{q})|$. Orient the discs so that there is one positive
trajectory between $\tilde{p}$ and $\tilde{q}$. Hence $[p:q]_{v'}=1-a$ with $\|a\|<1$. We get the term $a$
because there might be trajectories between $\tilde{p}$ and $g\tilde{q}$ with $\chi(g)<0$.\\
Conditions (2), (3) and (4) follow because the stable and unstable manifolds $W^{s,u}$ $(\tilde{r},\tilde{v}'')$
did not change in $f^{-1}([f(\tilde{r})-A_\mathcal{B}+L,f(\tilde{r})+A_\mathcal{B}-L])$ for critical points
$r\not=p$ or $q$ and since $v'$ is close enough to $v''$.
\end{proof}
\begin{theo}\label{cancalg}
Let $\omega$ be a Morse form without critical points of index $0,1,n-1,n$ which is $CC^1$. Let $v$ be a
transverse $\omega$-gradient and $L<0$. Let $p,\,q$ be critical points with ${\rm ind}\,p=
{\rm ind}\,q+1$ such that $[p:q]=\pm g(1-a)$ with $g\in G$ and $\|a\|<1$. Assume that $n=\dim M\geq 6$.
Then there is a Morse form
$\omega'$ cohomologous to $\omega$ with ${\rm crit}\,\omega'={\rm crit}\,\omega-\{p,q\}$ and a transverse
$\omega'$-gradient $v'$ such that:
\begin{enumerate}
\item $p_L([p':q']_{v'})=p_L([p':q']_v)$ for all $p',\,q'\in {\rm crit}\,\omega'$ with ${\rm ind}\,
p'={\rm ind}\,q'+1\not=i+1$.
\item $p_L(\psi_{v',v}(p'))=p'$ for all $p'\in {\rm crit}\,\omega'$.
\end{enumerate}
\end{theo}
\begin{proof}
Let $\tilde{p}$ and $\tilde{q}$ be the lifts of $p$ and $q$ used for the Novikov complex. By replacing
$\tilde{q}$ by $g\tilde{q}$ we can assume $[p:q]=\pm 1-a$. There exist only finitely many trajectories
between $\tilde{p}$ and all $h\tilde{q}$ with $\chi(h)\geq 0$. Furthermore, for every trajectory $T_1$
between $\tilde{p}$ and $\tilde{hq}$ with $\chi(h)\geq 0$ and $h\not=1$ there is another trajectory $T_2$
between $\tilde{p}$ and $\tilde{hq}$ with $\varepsilon(T_2)=-\varepsilon(T_1)$ since $[p:q](h)=0$.
So we can cancel such trajectories using Proposition \ref{canceltraj} provided there are no trajectories
between $\tilde{p}$ and $h'\tilde{q}$ with $\chi(h')>\chi(h)$. Start with the biggest $\chi(h)$ and cancel
all trajectories between $\tilde{p}$ and $h\tilde{q}$ with $\chi(h)\geq 0$ and $h\not=1$. Since $[p:q](1)=
\pm 1$ we can cancel all trajectories between $\tilde{p}$ and $\tilde{q}$ except one. Now we can cancel
$p$ and $q$ using Lemma \ref{cancfor}. The new transverse $\omega'$-gradient $v'$ can be arbitrarily
close to $v$ outside a small neighborhood of the trajectory so we can achieve conditions (1) and (2) since
the stable manifolds of critical points of index $\leq i$ with respect to $v$ stay away from the trajectory
and so do the unstable manifolds of critical points of index $\geq i+1$.
\end{proof}
\begin{rema}\label{remarko}\em
We do not obtain condition (1) of Theorem \ref{cancalg} in dimension $i+1$ as there can be trajectories of
$-v$ from a critical point $p'$ of index $i+1$ to $q$. After cancelling $q$ with $p$ these trajectories flow
towards $p$ under $-v'$ and from there to other critical points of index $i$ which appear in the boundary
of $p$ under $v$.
\end{rema}
\section{The simple homotopy type of the Novikov complex}\label{sec6}
In this section we assume that $M$ is a closed connected smooth manifold with $n=\dim M\geq 6$ and $\alpha
\in H^1(M;\mathbb{R})$ is $CC^1$. Given a finitely generated free $\widehat{\mathbb{Z}G}_\chi$ complex
$D_\ast$ with $D_i=0$ for $i\leq 1$ and $i\geq n-1$ which is simple chain homotopy equivalent to
$C^\Delta_\ast(M;\widehat{\mathbb{Z}G}_\chi)$, we want to realize it as the
Novikov complex of a Morse form representing $\alpha$. We will not be quite able to do this, but we can
approximate this in a reasonable sense. Notice that such complexes exist by Proposition \ref{no1n-1}.\\
If $A$ is a matrix over $\widehat{\mathbb{Z}G}_\chi$, denote $\|A\|=\max\{\|A_{ij}\|\}$ to be the norm of
$A$. This norm has similar properties as the norm for elements of $\widehat{\mathbb{Z}G}_\chi$, in particular
we have $\|A+B\|\leq \max\{\|A\|,\|B\|\}$ and $\|AB\|\leq \|A\|\,\|B\|$.
\begin{defi}\em
An invertible matrix $A$ over $\widehat{\mathbb{Z}G}_\chi$ is called \em simple \em if $\tau(A)=0\in
{\rm Wh}(G;\!\chi)$.
\end{defi}
Being a simple matrix is an open condition in the following sense:
\begin{lem}\label{opensim}
If $A$ is invertible, there is an $R_A>0$ such that $A-B$ is also invertible and $\tau(A-B)=\tau(A)$ for
$\|B\|<R_A$. If $A$ is simple, then $A=X(I-D)\in GL(\widehat{\mathbb{Z}G}_\chi)$ with $\|D\|\leq 1$ and
$\tau(X)=0\in K_1(\widehat{\mathbb{Z}G}_\chi)/\langle\pm[g]\,|\,g\in G\rangle$.
\end{lem}
\begin{proof}
We have $A-B=A(I-A^{-1}B)$, so choosing $R_A=\|A^{-1}\|^{-1}$ gives the first part. If $A$ is simple,
we have $\left(\begin{array}{cc}A&0\\0&I\end{array}\right)=E_1\cdots E_k$ with $E_i$
either an elementary matrix or a stabilization of $\pm g$ or $1-c$ with $\|c\|<1$. We can move matrices of
the form $1-c$ to the right of the product to get $\left(\begin{array}{cc}A&0\\0&I\end{array}\right)=X\cdot
(I-D)$ with $\|D\|<1$ and $X$ a product of elementary matrices and stabilizations of $\pm g$.
\end{proof}
Given a chain map $\varphi:D_\ast\to E_\ast$ between finitely generated free $\widehat{\mathbb{Z}G}_\chi$
complexes with given basis we can express each $\varphi_i:D_i\to E_i$ by a matrix which we also denote
by $\varphi_i$. Finitely generated free $\widehat{\mathbb{Z}G}_\chi$ complexes are assumed to have a basis,
in case of a Novikov complex the basis comes from liftings of critical points.
\begin{theo}\label{mainth}
Let $M^n$ be a closed connected smooth manifold with $n\geq 6$
and let $\alpha\in H^1(M;\mathbb{R})$ be
$CC^1$. Let $D_\ast$ be a finitely generated free $\widehat{\mathbb{Z}G}_\chi$ complex with $D_i=0$ for
$i\leq 1$ and $i\geq n-1$ which is simple chain homotopy equivalent to $C^\Delta_\ast(M;
\widehat{\mathbb{Z}G}_\chi)$. Given $L<0$ there is a Morse form $\omega$ representing
$\alpha$, a transverse $\omega$-gradient $v$ and a simple chain isomorphism $\varphi:D_\ast\to
C_\ast(\omega,v)$ where each $\varphi_i$ is of the form $I-A_i$ with $\|A_i\|<\exp L$.
\end{theo}
The condition that $\alpha$ be $CC^1$ cannot be removed as is shown in Damian \cite{damian}, where for
every $n\geq 8$ a manifold $M$ of dimension $n$ and a cohomology class $\alpha\in H^1(M;\mathbb{R})$ is
constructed such that the Novikov complex is simple chain homotopy equivalent to the trivial complex, but
$\alpha$ cannot be realized by a nonsingular 1-form.\\
Let us first outline the idea of the proof.
Given $D_\ast$, we choose any Novikov complex corresponding to a closed 1-form $\omega$. Then we introduce
for every generator of $D_2$ a pair of critical points of index 2 and 3. The new critical points of index
2 do not carry any information, but we can change the chain equivalence between $D_\ast$ and the new
Novikov complex so that the part between $D_2$ and the new critical points of index 2 approximates the
identity. Then we change the stable discs of these new critical points so that they carry the information
of the old critical points. Then the old critical points do not carry any relevant information and can be
traded against critical points of index 4. This way we can work our way up inductively until the Novikov
complex looks like $D_\ast$ except in dimensions $n-3$ and $n-2$.\\
Now we introduce for every generator
of $D_{n-3}$ and $D_{n-2}$ pairs of critical points of index $n-3$ and $n-2$, one which will carry the
information of the complex $D_\ast$ and one which is useless. Then we are left with the old critical points
of index $n-3$ and $n-2$ and the new useless critical points. The fact that $D_\ast$ has the simple
homotopy type of the Novikov complex now allows us to cancel these unnecessary critical points and we are
left with a Novikov complex which approximates $D_\ast$. In fact the boundary between the unnecessary
critical points forms a simple matrix which can be transformed to a matrix of the form $I-B$ with $\|B\|<1$.
by elementary steps. But this is good enough to cancel these critical points.\\
Before we start with the proof we need two algebraic lemmas first.
\begin{lem}\label{chan2iso}
Let $D_\ast,\,E_\ast$ be chain complexes, $\varphi:D_\ast\to E_\ast$ a chain map and $j$ an integer.
Assume that $E_j=C_j\oplus D_j$ and $E_{j+1}=C_{j+1}\oplus D_j$. Denote
\[\partial_{j+1}^E=\left(\begin{array}{cc}\partial_{11}&\partial_{12}\\ \partial_{21}&\partial_{22}
\end{array}\right)\hspace{1cm}\varphi_j=\left(\begin{array}{c}A_1\\A_2\end{array}\right)\hspace{1cm}
\varphi_{j+1}=\left(\begin{array}{c}B_1\\B_2\end{array}\right)\]
and assume that $\partial_{22}:D_j\to D_j$ is invertible. Define $\psi:D_\ast\to E_\ast$ by $\psi_i=\varphi_i$
for $i\not=j,\,j+1$ and
\[\psi_j=\left(\begin{array}{c}A_1+\partial_{12}\partial_{22}^{-1}\\ I+A_2\end{array}\right)\hspace{1cm}
\psi_{j+1}=\left(\begin{array}{c}B_1\\ \partial_{22}^{-1}\partial^D_{j+1}+B_2\end{array}\right)\]
where $I:D_j\to D_j$ denotes the identity. Then $\psi$ is chain homotopic to $\varphi$.
\end{lem}
\begin{proof}
Define $H:D_\ast\to E_{\ast+1}$ by $H_i=0$ for $i\not=j$ and $H_j=\left(\begin{array}{c}0\\
-\partial_{22}^{-1}\end{array}\right)$. Then
\[\partial^E_{j+1}H_j+H_{j-1}\partial^D_j=\partial^E_{j+1}H_j=\left(\begin{array}{cc}-\partial_{12}
\partial_{22}^{-1}\\ -I\end{array}\right)=\varphi_j-\psi_j\]
and
\[\partial^E_{j+2}H_{j+1}+H_j\partial^D_{j+1}=H_j\partial^D_{j+1}=\left(\begin{array}{c}0\\
-\partial_{22}^{-1}\partial^D_{j+1}\end{array}\right)=\varphi_{j+1}-\psi_{j+1}.\]
Hence $H$ is the required chain homotopy.
\end{proof}
\begin{lem}\label{isoinv}
Let $D_\ast,\,E_\ast$ be chain complexes, $j$ an integer and $\varphi:D_\ast\to E_\ast$ a chain homotopy
equivalence such that $\varphi_i:D_i\to E_i$ is an isomorphism for $i\leq j-1$. Then there is an
inverse equivalence $\psi:E_\ast\to D_\ast$ such that $\psi_i=\varphi_i^{-1}$ for $i\leq j-1$.
\end{lem}
\begin{proof}
Let $\psi':E_\ast\to D_\ast$ be a chain equivalence with ${\rm id}\,\simeq\psi'\varphi$ and ${\rm id}\,
\simeq\varphi\psi'$. Let $H:D_\ast\to D_{\ast+1}$ be a chain homotopy $H:{\rm id}\,\simeq \psi'\varphi$.
Define $\psi:E_\ast\to D_\ast$ by $\psi_i=\varphi^{-1}_i$ for $i\leq j-1$, $\psi_j=\psi_j'+H_{j-1}
\varphi^{-1}_{j-1}\partial^E_j$ and $\psi_i=\psi'_i$ for $i\geq j+1$. Now define $K:E_\ast\to D_{\ast+1}$
by $K_i=H_i\varphi_i^{-1}$ for $i\leq j-1$ and $K_i=0$ for $i\geq j$. Then it is easy to see that
$\partial^D_{i+1}K_i+K_{i-1}\partial^E_i=\psi_i-\psi_i'$. for all $i$.
\end{proof}
\begin{proof}[Proof of Theorem \ref{mainth}]
By Proposition \ref{no1n-1} there is a Morse form $\omega$ representing $\alpha$ without critical points of
index $0,1,n-1,n$. Choose any transverse $\omega$-gradient $v$. The Novikov complex $C_\ast(\omega,v)$
is simple chain homotopy equivalent to $C_\ast^\Delta(M;\widehat{\mathbb{Z}G}_\chi)$,
so there is a simple chain homotopy equivalence $\varphi:D_\ast\to C_\ast(\omega,v)$. Denote
$C_\ast=C_\ast(\omega,v)$.\\
Assume we have $j\leq n-4$ such that we have a simple chain homotopy equivalence $\varphi:D_\ast\to C_\ast$
such that $\varphi_i=I-A_i$ with $\|A_i\|<\exp L$ for $i\leq j-1$. Note that this is true for $j=2$. We
want to find a new Morse form such that this is also true for $j+1$.\\
{\bf Step 1: Introduction of new critical points of index \boldmath$j$}\qua
Let $q^j_1,\ldots,q_{k_j}^j$ be the critical points of $\omega$ with index $j$ and let $d^j_1,\ldots
,d^j_{l_j}$ be the generators of $D_j$. Denote the chain inverse of $\varphi$ by $\bar{\varphi}$. By Lemma
\ref{isoinv} we can assume that $\bar{\varphi}_i=\varphi_i^{-1}$ for $i\leq j-1$. Also let $H:C_\ast\to
C_{\ast+1}$ be a chain homotopy $H:{\rm id}\,\simeq \varphi\circ\bar{\varphi}$. Since $C_\ast$ is free,
we can assume that $H_i=0$ for $i\leq j-1$, compare Dold \cite[Ex.VI.1.12.4]{dold}.\\
For every $d_l^j$ we introduce a pair of critical points $p_l^j$
and $p^{j+1}_l$ of index $j$ and $j+1$ by Proposition \ref{insert}, thus getting a new Morse form
$\omega'$ and a transverse $\omega'$-gradient $v'$. Also, we can achieve this so that $(\psi_{v,v'})_j=
\left(\begin{array}{c}I-E_j\\E_j'\end{array}\right)$, $(\psi_{v,v'})_{j+1}=\left(\begin{array}{c}I-E_{j+1}\\
E_{j+1}'\end{array}\right)$ and $(\psi_{v,v'})_i=I-E_i$ for $i\not=j,\,j+1$ with $\|E_i\|,\,\|E_i'\|$ $<
\min\{1,\|\varphi\|^{-1}\}\cdot\exp L$. Also $\partial'_{j+1}:C_{j+1}(\omega',v')\to C_j(\omega',v')$ is of
the form
\begin{eqnarray*}
\partial'_{j+1}&=&\left(\begin{array}{cc}\partial_{j+1}-F_1&F_2\\F_3&I-A\end{array}\right)
\end{eqnarray*}
with $\partial_{j+1}:C_{j+1}(\omega,v)\to C_j(\omega,v)$ the boundary, $\|A\|<1$ and $\|F_i\|<\min\{1,\|
\varphi\|\}\cdot
\min\{1,\|H\|^{-1}, \|\bar{\varphi}\|^{-1}\}$. The composition of $\varphi$ and $\psi_{v,v'}$
gives a simple chain homotopy
equivalence $\varphi'$ with $\|\varphi_i'-\varphi_i\|<\exp L$ for $i\not=j,\,j+1$ and $\varphi_j'=\left(
\begin{array}{c}\varphi_j-F_j\\F_j'\end{array}\right)$, $\varphi_{j+1}'=\left(\begin{array}{c}\varphi_{j+1}
-F_{j+1}\\F_{j+1}'\end{array}\right)$, where $\|F_i\|,\,\|F_i'\|<\exp L$. So by Lemma \ref{chan2iso} we
have a simple chain homotopy equivalence $\psi:D_\ast\to C_\ast(\omega',v')$ with $\psi_j=\left(
\begin{array}{c}\varphi_j-F_j''\\I+F_j'\end{array}\right)$, $\psi_{j+1}=\left(\begin{array}{c}\varphi_{j+1}
-F_{j+1}\\X\end{array}\right)$ with $\|F_j''\|<\exp L$ and $\psi_i=\varphi'_i$ for $i\not=j,\,j+1$.
Since $L<0$, $I+F_j'$ is invertible.\\
Let $E_\ast=C_\ast(\omega',v')$. Then $E_j=C_j\oplus D_j$. Perform an elementary change of basis on $E_j$
of the form $\left(\begin{array}{cc}I&-(\varphi_j-F_j'')(I+F_j')^{-1}\\0 & I\end{array}\right)$. With this
change of basis the matrix of $\psi_j$ is of the form $\left(\begin{array}{c}0\\I+F_j'\end{array}\right)$.
Using Proposition \ref{elmchan}, we can approximate this elementary change of basis arbitrary well.
So approximate the elementary change of basis so that we get a Morse form $\omega''$, a transverse
$\omega''$-gradient $v''$ and a simple chain homotopy equivalence $\psi':D_\ast\to C_\ast(\omega'',v'')$ with
$\psi_j'=\left(\begin{array}{c}Y\\I-G_j\end{array}\right)$, $\psi_{j+1}'=\left(\begin{array}{c}
\varphi_{j+1}-G_{j+1}\\Z\end{array}\right)$, and $\|\psi'_i-\varphi_i\|<\exp L$ for $i\not=j,\,j+1$ with
$\|G_i\|<\exp L$ and the boundary $\partial''_{j+1}:C_{j+1}(\omega'',v'')\to C_j(\omega'',v'')$ is of
the form
\begin{eqnarray*}
\partial''_{j+1}&=&\left(\begin{array}{cc}\partial_{j+1}-K_1&-\varphi_j-K_2\\-K_3&I-A'\end{array}\right)
\end{eqnarray*}
with $\|K_i\|<\min\{1,\|H\|^{-1},\|\bar{\varphi}\|^{-1}\}$ and $\|A'\|<1$.\\
Define $L_j:C_j(\omega'',v'')\to C_{j+1}(\omega'',v'')$ by $L_j=\left(\begin{array}{cc}H_j&0\\
-\bar{\varphi}_j&0\end{array}\right)$. Then
\begin{eqnarray*}
\partial''_{j+1}\circ L_j&=&\left(\begin{array}{cc}
\partial_{j+1}H_j-K_1H+\varphi_j\bar{\varphi}_j+K_2\bar{\varphi}_j&0\\-K_3H_j-\bar{\varphi}_j-A'
\bar{\varphi}_j&0\end{array}\right)\\
&=&\left(\begin{array}{cc}I-S&0\\ \ast&0\end{array}\right)
\end{eqnarray*}
with $\|S\|<1$. So for every critical point $q_k^j$ there exists a $u_k\in C_{j+1}(\omega'',v'')$ with
$\partial''_{j+1}u_k=(1-a_k)q^j_k+r$ with $\|a_k\|<1$ and $p_{q,k}(r)=0$, where $p_{q,k}:C_j(\omega'',v'')\to
C_\omega'',v'')$ is projection to the span of $q^j_k$.\\
Rename $\omega=\omega''$, $v=v''$ and $\varphi=\psi'$.\\
{\bf Step 2: Removal of unnecessary critical points of index \boldmath$j$}\qua
The critical points of $\omega$ of index $j$ are $q^j_1,\ldots,q^j_{k_j}$ and $p_1^j,\ldots,p^j_{l_j}$ where
the $q^j_k$ are the critical points of the original $\omega$ and the $p^j_k$ correspond to the generators
$d^j_k$ of $D_j$. For $q_k^j$ introduce a pair of critical points $r^{j+1}_k$ and $r^{j+2}_k$ of index
$j+1$ and $j+2$ with Proposition \ref{insert} to get a new Morse form $\omega'$ and a transverse
$\omega'$-gradient $v'$ so that $\varphi'=\psi_{v,v'}\circ\varphi$ satisfies $\varphi'_i=\varphi_i-E_i$
for $i\not=j+1,\,j+2$, $\varphi_{j+1}'=\left(\begin{array}{c}\varphi_{j+1}\\E_{j+1}\end{array}\right)$ and
$\varphi_{j+2}'=\left(\begin{array}{c}\varphi_{j+2}\\E_{j+2}\end{array}\right)$ with $\|E_i\|<\exp L$ for
all $i$. In particular we have $p_{q,k}(\partial r^{j+1}_k)=a q^j_k$ with $\|a\|<\exp L$.\\
With the elementary change of basis on $C_{j+1}(\omega',v')$ of the form $\bar{r}^{j+1}_k=r^{j+1}_k+u_k$
we get $p_{q,k}(\partial \bar{r}^{j+1}_k)=(1-a_k+a)q^j_k$. So use Proposition \ref{elmchan} to get a new
Morse form $\omega''$ and transverse $\omega''$-gradient $v''$ such that for the critical point $r_k^{j+1}$
we now have $p_{q,k}(\partial r_k^{j+1})=(1-b)q^j_k$ with $\|b\|<1$. Therefore we can cancel the critical
points $r_k^{j+1}$ and $q^j_k$ for all $k$ using Theorem \ref{cancalg}. Remember we have $\varphi_j=
\left(\begin{array}{c}Y\\I-G_j\end{array}\right)$ with $\|G_j\|<\exp L$. We can cancel so that for the new
Morse form without the critical points $q^j_1,\ldots, q^j_{k_j}$ we now have $\varphi''_i=I-G_i'$ with
$\|G_i'\|<\exp L$ for all $i\leq j$. Therefore we have finished the induction step.\\
So we can assume that we have a simple chain homotopy equivalence $\varphi:D_\ast\to C_\ast(\omega,v)$
such that $\varphi_i=I-A_i$ with $\|A_i\|<\exp L$ for $i\leq n-4$. Notice also that everything we have done
so far would have worked if $D_\ast$ was just chain homotopy equivalent to the Novikov complex. But to get
the result in the final two dimensions, we need the same simple homotopy type. Denote $C_\ast=C_\ast(\omega,
v)$.\\
{\bf Step 3: Introduction of new critical points in dimension \boldmath$n-3$ and $n-2$}\qua
We want to introduce new critical points of index $n-3$ and $n-2$ for every generator of $D_{n-3}$ and
$D_{n-2}$. Let us do this on an algebraic level first. We have a simple chain homotopy equivalence
$\varphi:D_\ast\to C_\ast$ such that $\varphi_i:D_i\to C_i$ is a simple isomorphism for $i\leq n-4$.
Define a new chain complex $E_\ast$ by $E_i=C_i$ and $\partial_i^E=\partial^C_i$ for $i\leq n-4$,
$E_{n-3}=C_{n-3}\oplus D_{n-2}\oplus D_{n-3}$, $E_{n-2}=C_{n-2}\oplus D_{n-3}\oplus D_{n-2}$ and
\[\partial^E_{n-2}=\left(\begin{array}{ccc}\partial^C_{n-2}&0&0\\0&0&I-A_{n-2}\\0&I-A_{n-3}&0\end{array}
\right)\hspace{1cm}\partial^E_{n-3}=\left(\begin{array}{ccc}\partial^C_{n-3}&0&0\end{array}\right)\]
with $\|A_{n-2}\|,\,\|A_{n-3}\|<1$. It is easy to see that $E_\ast$ is simple homotopy equivalent to $C_\ast$
and $\psi:D_\ast\to E_\ast$ defined by $\psi_i=\varphi_i$ for $i\leq n-4$, $\psi_{n-3}=\left(\begin{array}{c}
\varphi_{n-3}\\0\\0\end{array}\right)$, $\psi_{n-2}=\left(\begin{array}{c}\varphi_{n-2}\\0\\0\end{array}
\right)$ is a simple chain homotopy equivalence.\\
By Lemma \ref{chan2iso} $\psi$ is chain homotopic to $\psi'$ with $\psi'_i=\psi_i$ for $i\leq n-4$,
$\psi'_{n-3}=\left(\begin{array}{c}\varphi_{n-3}\\0\\I\end{array}\right)$, $\psi'_{n-2}=\left(
\begin{array}{c}\varphi_{n-2}\\(I-A_{n-3})^{-1}\partial^D_{n-2}\\0\end{array}\right)$. Let
$\bar{\varphi}$ be a chain inverse to $\varphi$ such that $\bar{\varphi}_i=\varphi^{-1}$ for $i\leq n-4$
and $K:D\ast\to D_{\ast+1}$ a chain homotopy $K:{\rm id}\,\simeq\bar{\varphi}\circ\varphi$ such that
$K_i=0$ for $i\not= n-3$.\\
Now define a chain homotopy $H_\ast:D_\ast\to E_{\ast+1}$ by $H_i=0$ for $i\not=n-3$ and $H_{n-3}=\left(
\begin{array}{c}0\\0\\K_{n-3}\end{array}\right)$. Then $H:\psi'\simeq\psi''$ with $\psi''_i=\psi'_i$ for
$i\leq n-4$ and \\ 
$\psi_{n-3}''=\left(\begin{array}{c}\varphi_{n-3}\\(I-A_{n-2})K_{n-3}\\I\end{array}\right)$,
$\psi''_{n-2}=\left(\begin{array}{c}\varphi_{n-2}\\(I-A_{n-3})^{-1}\partial^D_{n-2}\\I-\bar{\varphi}_{n-2}
\varphi_{n-2}\end{array}\right)$.
Perform a change of basis on $E_{n-2}$ of the form $\left(\begin{array}{ccc}I&0&0\\0&I&0\\
\bar{\varphi}_{n-2}&0&I\end{array}\right)$, then the matrix of $\psi''_{n-2}$ is $\left(\begin{array}{c}
\varphi_{n-2}\\(I-A_{n-3})^{-1}\partial^D_{n-2}\\I\end{array}\right)$ and the boundary matrix is
\[\partial^E_{n-2}=\left(\begin{array}{ccc}\partial^C_{n-2}&0&0\\
-(I-A_{n-2})\bar{\varphi}_{n-2}&0&I-A_{n-2}\\
0&I-A_{n-3}&0\end{array}\right)\]
Define $F_i=0$ for $i\leq n-4$, $F_{n-3}=C_{n-3}\oplus D_{n-2}$ and $F_{n-2}:C_{n-2}\oplus D_{n-3}$. Now
$\psi''$ is a chain map $\psi''=\left(\begin{array}{c}\psi''_F\\\psi''_D\end{array}\right):D_i\to F_i\oplus
D_i$ with $\psi''_D$ a simple automorphism for every $i$. By Ranicki \cite[Prop.1.8]{ranick} we have that
${\rm coker}(\psi'')$ is isomorphic to a chain complex $\hat{F}_\ast$ with $\hat{F}_i=F_i$ and
\begin{eqnarray*}
\partial^{\hat{F}}_{n-2}&\!=\!&\left(\begin{array}{cc}\partial_{n-2}^C&0\\-(I-A_{n-2})\bar{\varphi}_{n-2}&0
\end{array}\right)-\left(\begin{array}{c}\varphi_{n-3}\\(I-A_{n-2})K_{n-3}\end{array}\right)\left(
\begin{array}{cc}0&I-A_{n-3}\end{array}\right)\\
&\!=&\!\left(\begin{array}{cc}I&0\\0&I-A_{n-2}\end{array}\right)\left(\begin{array}{cc}\partial^C_{n-2}
&-\varphi_{n-3}\\
-\bar{\varphi}_{n-2}&-K_{n-3}\end{array}\right)\left(\begin{array}{cc}I&0\\0&I-A_{n-3}\end{array}\right)
\end{eqnarray*}
Furthermore, again by Ranicki \cite[Prop.1.8]{ranick} the natural projection $p:\mathcal{C}(\psi'')\to {\rm coker}(\psi'')=\hat{F}$ is a chain homotopy equivalence with torsion
\[\tau(p)=\sum_{i=2}^{n-2}(-1)^{i+1}\tau(\psi_D'':D_i\to D_i)\in {\rm Wh}(G;\chi),\]
so $\tau(p)=0$ and we have $\tau(\hat{F})=
\tau(\psi'')=0$. Therefore 
\[D:=\left(\begin{array}{cc}\partial^C_{n-2}&-\varphi_{n-3}\\-\bar{\varphi}_{n-2}&
-K_{n-3}\end{array}\right)\] 
is a simple matrix. By Lemma \ref{opensim} there is an $R>0$ such that $D-B$
is also simple for $\|B\|<R$. Also $\partial^{\hat{F}}_{n-2}-B$ is simple for $\|B\|<R$ and $\|A_{n-2}\|,
\,\|A_{n-3}\|<1$.
Now perform a change of basis on $E_{n-3}$ of the form $\left(\begin{array}{ccc}I&0&-\varphi_{n-3}\\0&I&
-(I-A_{n-2})K_{n-3}\\0&0&I\end{array}\right)$. Then the matrix of $\partial^E_{n-2}$ is
\[\partial^E_{n-2}=\!\left(\!\begin{array}{ccc}\partial^C_{n-2}&-\varphi_{n-3}(I-A_{n-3})&0\\-(I-A_{n-2})
\bar{\varphi}_{n-2}&-(I-A_{n-2})K_{n-3}(I-A_{n-3})&I-A_{n-2}\\0&I-A_{n-3}&0\end{array}\!\right)\!=:\bar{D}\]
Now introduce as in Step 1 new critical points $p^{n-3}_t$ and $r_t^{n-2}$ of index $n-3$ and $n-2$ for
every generator of $D_{n-3}$ and critical points $p^{n-2}_k$ and $r^{n-3}_k$ of index $n-2$ and $n-3$
for every generator of $D_{n-2}$ to get a new Morse form $\omega'$ and transverse $\omega'$-gradient
$v'$. We can approximate the described change of basis on $C_\ast(\omega',v')$ to end up with a Morse
form $\omega''$ and a Novikov complex $C_\ast(\omega'',v'')$ such that $\partial_{n-2}=\bar{D}-X$ with
$\|X\|$ arbitrary small. In particular we can make it so small that the submatrix, denoted $\bar{\partial}$,
corresponding to the critical points $\{q_s^{n-2},\,r^{n-2}_t\}$ and $\{q_l^{n-3},\,r_k^{n-3}\}$ is simple
and $\psi_i=I-A_i$, $\psi_{n-3}=\left(\begin{array}{c}\ast\\ \ast\\ I-A_{n-3}\end{array}\right)$ and
$\psi_{n-2}=\left(\begin{array}{c}\ast\\ \ast \\I-A_{n-2}\end{array}\right)$ with $\|A_i\|<\exp L$ for
all $i$.\\
{\bf Step 4: Elimination of critical points in dimension \boldmath$n-3$ and $n-2$}\qua
Using Lemma \ref{opensim} we can change $\bar{\partial}$ into a matrix of the form $I-B$ with $\|B\|<1$ by
elementary changes of basis and stabilizing. Approximate these changes of basis and add critical points
so that for the Novikov complex we have
\begin{eqnarray*}\partial_{n-2}&=&\left(\begin{array}{cc}I-B'&\ast\\ \ast& \ast\end{array}\right)
\end{eqnarray*}
with $\|B'\|<1$. Now we can cancel all critical points $\{q_l^{n-3},\,r^{n-3}_k\}$ against the critical
points $\{q_s^{n-2},\,r_t^{n-2}\}$ to get the required Morse form.
\end{proof}
As a corollary we get Latour's theorem \cite{latour}. If the chain complex $C^\Delta_\ast(M;
\widehat{\mathbb{Z}G}_\chi)$ is acyclic, define the \em Latour obstruction \em to be
$\tau(M,\alpha)=\tau(C^\Delta_\ast(M;\widehat{\mathbb{Z}G}_\chi))\in{\rm Wh}(G;
\chi)$.
\begin{theo}\label{latthe}
Let $M^n$ be a closed connected smooth manifold with $n\geq 6$ and $\alpha\in H^1(M;\mathbb{R})$.
Then $\alpha$ can be represented by a closed 1-form without critical points if and only if $\alpha$ is
$CC^1$, $C^\Delta_\ast(M;\widehat{\mathbb{Z}G}_\chi)$ is acyclic and $\tau(
M,\alpha)=0\in {\rm Wh}(G;\chi)$.
\end{theo}
\begin{rema}\em
To proof Theorem \ref{latthe} directly, notice that the fairly involved steps 1 and 3 are not needed for
this.
\end{rema}
To compare Theorem \ref{mainth} to Pajitnov \cite{pajisu} we need a new notion.
\begin{defi}\em
Let $N\in\mathbb{R}$. Two finitely generated free $\widehat{\mathbb{Z}G}_\chi$ chain complexes with basis
and ${\rm rank}\, D_i={\rm rank}\, E_i$ for all $i$ are called \em $N$-equivalent \em if $\|\partial^D-
\partial^E\|\leq
\exp N$.
\end{defi}
The analogue of Pajitnov \cite[Th.0.12]{pajisu} is as follows.
\begin{theo}\label{pajith}
Let $M^n$ be a closed connected smooth manifold with $n\geq 6$ and $\alpha\in H^1(M;\mathbb{R})$ be
$CC^1$. Let $D_\ast$ be a finitely generated free $\widehat{\mathbb{Z}G}_\chi$ complex with $D_i=0$ for
$i\leq 1$ and $i\geq n-1$ which is simple chain homotopy equivalent to $C^\Delta_\ast(M;
\widehat{\mathbb{Z}G}_\chi)$. Given $N\in\mathbb{R}$ there is a Morse form $\omega$
representing $\alpha$ and a transverse $\omega$-gradient $v$ such that $D_\ast$ and $C_\ast(\omega,v)$
are $N$-equivalent.
\end{theo}
\proof
By Theorem \ref{mainth} there is a Morse form $\omega$, a transverse $\omega$-gradient $v$ and a simple
chain isomorphism $\varphi:D_\ast\to C_\ast(\omega,v)$ with $\varphi_i=I-A_i$ and $\|A_i\|<
\|\partial\|^{-1}\exp N$. Then, since $\varphi$ is a chain map, we have $\varphi_{i-1}\partial_i^D
\varphi_i^{-1}=\partial^C_i$. As matrices we get
\begin{align*}
\|\partial^D_i-\partial^C_i\|&=\|\partial^D_i-(I-A_{i-1})\partial^D_i(I-A_i)^{-1}\|\\
&\leq\max\{\|A_{i-1}\partial^D_i\|,\,\|\partial^D_iA_i\|\}\\
&\leq\|\partial^D_i\|\max\{\|A_{i-1}\|,\,\|A_i\|\}\\
&\leq\|\partial^D_i\|\cdot\|\partial^D_i\|^{-1}\cdot\exp N.\tag*{\qed}
\end{align*}

Instead of the Novikov ring we can look at a certain noncommutative Cohn localization. Let $\Sigma_\chi$
be the set of diagonal matrices over $\mathbb{Z}G$ of the form $I-A$ with $\|A\|_\chi<1$. By Cohn
\cite{cohn} there is a unique ring $\Sigma^{-1}_\chi\mathbb{Z}G$ and a natural ring homomorphism $i:
\mathbb{Z}G\to\Sigma^{-1}_\chi\mathbb{Z}G$ with $i(\Sigma_\chi)\subset GL(\Sigma^{-1}_\chi\mathbb{Z}G)$
such that for every ring homomorphism $\eta:\mathbb{Z}G\to R$ with $\eta(\Sigma_\chi)\subset GL(R)$
there is a unique ring homomorphism $\varepsilon:\Sigma^{-1}_\chi\mathbb{Z}G\to R$ with $\varepsilon\circ i
=\eta$.\\
Notice that the inclusion $j:\mathbb{Z}G\to\widehat{\mathbb{Z}G}_\chi$ satisfies $j(\Sigma_\chi)\subset
GL(\widehat{\mathbb{Z}G}_\chi)$, so there is a ring homomorphism $\varepsilon:\Sigma^{-1}_\chi\mathbb{Z}
G\to\widehat{\mathbb{Z}G}_\chi$ with $j=\varepsilon\circ i$. In particular $i:\mathbb{Z}G\to\Sigma^{-1}_\chi
\mathbb{Z}G$ is injective. Define ${\rm Wh}(G;\Sigma_\chi)=K_1(\Sigma^{-1}_\chi\mathbb{Z}G)/\langle[\pm g],
\,[i(\Sigma_\chi)]\rangle$.\\
A result of Farber \cite{farber} says that given a Morse form $\omega$ there is a finitely generated free
$\Sigma^{-1}_\chi\mathbb{Z}G$ complex $D_\ast$ simple chain homotopy equivalent to $C^\Delta_\ast(M;$
$\Sigma^{-1}_\chi\mathbb{Z}G)$ with ${\rm rank}\, D_i=c_i(\omega)=|\{p\in M\,|\,
\omega_p=0 \mbox{ and }{\rm ind}\,p=i\}|$. Notice that Farber \cite[Lm.8.12]{farber} points out that $D_\ast$
need not be simple chain homotopic to $C_\ast^\Delta(M;\Sigma^{-1}_\chi\mathbb{Z}G)$
when viewed over $K_1(\Sigma^{-1}_\chi\mathbb{Z}G)/\langle[\pm g]\rangle$. But by comparing the
proof of \cite[Lm.8.12]{farber} with \cite[Lm.7.1]{farber} and Ranicki \cite[Prop.1.8]{ranick} one sees
that the torsion of the last collapse in \cite[Lm.8.12]{farber} vanishes in ${\rm Wh}(G;\Sigma_\chi)$.
Combining this with Theorem \ref{mainth} we get the following.
\begin{theo}\label{localv}
Let $M$ be a closed connected smooth manifold with $n=\dim M\geq 6$ and $\alpha\in H^1(M;\mathbb{R})$ be
$CC^1$. Then:
\begin{enumerate}
\item Given a finitely generated free $\Sigma^{-1}_\chi\mathbb{Z}G$ complex $D_\ast$ with $D_i=0$ for
$i\leq 1$ and $i\geq n-1$ simple chain homotopy equivalent to $C^\Delta_\ast(M;
\Sigma^{-1}_\chi\mathbb{Z}G)$ there is a Morse form $\omega$ with $c_i(\omega)={\rm rank}\,
D_i$.
\item Given a finitely generated free $\widehat{\mathbb{Z}G}_\chi$ complex $E_\ast$ with $E_i=0$ for
$i\leq 1$ and $i\geq n-1$ simple chain homotopy equivalent to $C_\ast^\Delta(M;
\widehat{\mathbb{Z}G}_\chi)$ there is a finitely generated free $\Sigma^{-1}_\chi\mathbb{Z}G$ complex
$D_\ast$ with ${\rm rank}_{\Sigma_\chi^{-1}\mathbb{Z}G}D_i={\rm rank}_{\widehat{\mathbb{Z}G}_\chi}E_i$
simple chain homotopy equivalent to $C_\ast^\Delta(M;\Sigma^{-1}_\chi\mathbb{Z}G)$.
\end{enumerate}
\end{theo}
In particular the Latour obstruction for the existence of a closed 1-form without critical points pulls back
to an obstruction in ${\rm Wh}(G;\Sigma_\chi)$. In the rational case the obstruction actually pulls back to
${\rm Wh}(G)$, see the original fibering obstructions of Farrell \cite{farrel,farre2} or Siebenmann
\cite{sieben} and their comparison to the Latour obstruction in Ranicki \cite{ranict}. This raises the
question whether the Latour obstruction can be pulled back to an obstruction in ${\rm Wh}(G)$ in general.
\begin{rema}\em
Theorem \ref{mainth} reduces the problem of finding a Morse form with a minimal number of critical points
in a $CC^1$ cohomology class $\alpha$ on a manifold $M$ with dimension $\geq 6$ to the algebraic problem
of finding a finitely generated free $\widehat{\mathbb{Z}G}_\chi$ complex $D_\ast$ simple homotopy
equivalent to $C_\ast^\Delta(M;\widehat{\mathbb{Z}G}_\chi)$ with a minimal
number of generators and with $D_i=0$ for $i\leq 1$ and $i\geq n-1$. The last condition that $D_i=0$ for
$i\leq 1$ and $i\geq n-1$ can be removed using Pajitnov \cite[Prop.7.14]{pajisu}. By Theorem \ref{localv}
we can furthermore use $\Sigma^{-1}_\chi\mathbb{Z}G$ instead of $\widehat{\mathbb{Z}G}_\chi$.
\end{rema}
\section{Realization of torsion}
In this section we analyze the impact of Theorem \ref{mainth} on the
torsion of the chain homotopy equivalence $\varphi_v:C^\Delta_\ast(M;
\widehat{\mathbb{Z}G}_\chi)\to C_\ast(\omega,v)$ described in the
appendix. We know by Theorem \ref{simple} that the torsion vanishes in
${\rm Wh}(G;\chi)$, but it is known that $\tau(\varphi_v)$ is a well
defined element of the subgroup $\overline{W}$ of
$K_1(\widehat{\mathbb{Z}G}_\chi)/\langle[\pm g]\,|\,g\in G\rangle$ generated
by units of the form $1-a$ with $\|a\|<1$. This torsion also carries
information about the closed orbit structure of $v$ in form of a zeta
function, see \cite[Th.1.1]{schue2}. So realizing a given element of
$\overline{W}$ as the torsion of $\varphi_v$ for some combination of
$\omega$ and $v$ implies the realization of a zeta function. The
result we can prove now reads as follows.
\begin{theo}
Let $G$ be a finitely presented group, $\chi:G\to\mathbb{R}$ be
$CC^1$, $b\in\widehat{\mathbb{Z}G}_\chi$ satisfy $\|b\|<1$ and
$\varepsilon>0$. Then for any closed connected smooth manifold $M$
with $\pi_1(M)=G$ and $\dim M\geq 6$ there is a Morse form $\omega$
realizing $\chi$, a transverse $\omega$-gradient $v$ and a
$b'\in\widehat{\mathbb{Z}G}_\chi$ with $\|b-b'\|<\varepsilon$ such
that $\tau(\varphi_v)=\tau(1-b')\in
K_1(\widehat{\mathbb{Z}G}_\chi)/\langle[\pm g]\,|\,g\in G\rangle$.
\end{theo}
\begin{proof}
Choose a Morse form $\omega'$ representing $\chi$ and a transverse
$\omega'$-gradient $v'$. Let $1-c\in\widehat{\mathbb{Z}G}_\chi$
represent $\tau(1-b)-\tau(\varphi_{v'})\in\overline{W}$. Let $C_\ast
=C_\ast(\omega',v')$. Denote by $C(1-c)_\ast$ the finitely generated
free $\widehat{\mathbb{Z}G}_\chi$ complex with $C(1-c)_j=0$ for
$j\not=n-3,n-2$, where $n=\dim M$,
$C(1-c)_j=\widehat{\mathbb{Z}G}_\chi$ for $j=n-3,n-2$ and
$d:C(1-c)_{n-2}\to C(1-c)_{n-3}$ is multiplication by
$(1-c)^{(-1)^{n-1}}$. Then $C(1-c)_\ast$ is acyclic with
$\tau(C(1-c)_\ast)= \tau(1-b)-\tau(\varphi_{v'})$. Also $D_\ast=C_\ast
\oplus C(1-c)_\ast$ is simple homotopy equivalent to $C_\ast$. By
Theorem \ref{mainth} $D_\ast$ can be approximately realized as the
Novikov complex of a Morse form $\omega$ and a transverse
$\omega$-gradient $v$. Note that in the proof of Theorem \ref{mainth}
we can start directly with Step 3 and we only have to introduce
critical points for the generators of $C(1-c)_\ast$. By analyzing the
proof using Section \ref{sec5} we see that there is a sequence of
Morse forms $\omega_i$, $i=1,\ldots,k$ with $\omega_1=\omega'$,
$\omega_k=\omega$ and $\omega_i$ agrees with $\omega'$ in a
neighborhood of the critical points of $\omega'$. Furthermore there
are homotopy equivalences $\varphi^i:C_\ast(\omega_i,v_i)\to
C_\ast(\omega_{i+1},v_{i+1})$ chain homotopic to $\psi_{v_i,v_{i+1}}$
and the matrix of $\varphi^i$ restricted to the subgroup generated by
the critical points of $\omega'$ is of the form $I-A$ with $\|A\|<
\varepsilon$. Denote $\varphi=\varphi^{k-1}\circ\cdots\circ\varphi^1$,
then $\tau(\psi_{v',v})=\tau(\varphi)$ by Proposition
\ref{formcom}. We have $C_j(\omega,v)=C_j$ for $j\not=n-3,n-2$ and
$C_j(\omega,v)=C_j\oplus\widehat{\mathbb{Z}G}_\chi$ for
$j=n-3,n-2$. Since all $\varphi_j$ restricted to $C_j$ are of the form
$I-A_j$ with $\|A_j\|<\varepsilon$ we get that $\varphi_j$ is a split
injection and that $\mathcal{C}(\varphi)$ is chain homotopy equivalent
to ${\rm coker}(\varphi)$ by the projection $p:\mathcal{C}(\varphi)
\to {\rm coker}(\varphi)$, see Ranicki
\cite[Prop.1.8]{ranick}. Furthermore
$\tau(p)=\sum_{j=2}^{n-2}(-1)^{j+1}\tau(\varphi_j:C_j\to C_j)$. Also
${\rm coker}(\varphi)$ is an approximation of $C(1-c)_\ast$, i.e.\
$\tau({\rm coker}(\varphi))=\tau(1-c)+\tau(1-e)$ where
$e\in\widehat{\mathbb{Z}G}_\chi$ satisfies
$\|e\|<\varepsilon$. Therefore 
\[\tau(\psi_{v',v})=\tau(\varphi)=\tau({\rm
coker}(\varphi))-\tau(p)=\tau(1-c)-\tau(1-e')
\]
with $\|e'\|<\varepsilon$. By Proposition \ref{formcom} we now get
\[\tau(\varphi_v)=\tau(\psi_{v',v})+\tau(\varphi_{v'})=\tau(1-b)-\tau(\varphi_{v'})
-\tau(1-e')+\tau(\varphi_{v'})\]
This gives the result.
\end{proof}
\section{Poincar\'{e} duality}
Let $M$ be a closed connected smooth manifold, $\omega$ a Morse form and $v$ a transverse $\omega$-gradient.
Then $-\omega$ is a Morse form as well and $-v$ a transverse $(-\omega)$-gradient. To define the Novikov
complex $C_\ast(\omega,v)$ we need to choose orientations of $W^s(p,v)$ which induce coorientations of
$W^u(p,v)$ and liftings $\tilde{p}\in\tilde{M}$ for all critical points $p$ of $\omega$. These orientations
lift to orientations of $W^s(g\tilde{p},\tilde{v})$ for all $g\in G$. To define $C_\ast(-\omega,-v)$ we need
orientations for $W^s(p,-v)=W^u(p,v)$. The universal cover $\tilde{M}$ is orientable, so fix an orientation.
Denote chosen orientations by $o(N)$ for orientable manifolds $N$. Now choose for every critical point $p$
an orientation of $W^s(\tilde{p},-\tilde{v})$ such that $o(W^s(\tilde{p},\tilde{v}))\wedge o(W^s(\tilde{p},
-\tilde{v}))=o(\tilde{M})$, where the wedge means ''followed by''. Use the covering transformations to orient
$W^s(g\tilde{p},-\tilde{v})$ for all $g\in G$ and the projection to orient $W^s(p,-v)$. Then $o(W^s(
g\tilde{p},\tilde{v}))\wedge o(W^s(g\tilde{p},-\tilde{v}))=w(g)\cdot o(\tilde{M})$ where $w:G\to \{\pm 1\}$
is the orientation homomorphism of $M$.\\
Let $p,\,q$ be critical points of $\omega$ with ${\rm ind}\,p={\rm ind}\, q+1=i$. Let $T$ be a trajectory
between $\tilde{p}$ and $g\tilde{q}$, where $\tilde{p}$ and $\tilde{q}$ are the chosen liftings of $p$
and $q$. Then $g^{-1}(-T)$ is a trajectory between $\tilde{q}$ and $g^{-1}\tilde{p}$. With the choice of
orientations we now get
\begin{eqnarray}\label{orienter}
\varepsilon(g^{-1}(-T))&=&w(g)(-1)^i \varepsilon(T)
\end{eqnarray}
where $\varepsilon(T)$ and $\varepsilon(g^{-1}(-T))$ are defined as in Section \ref{sec5}.\\
The involution $\bar{ }:\mathbb{Z}G\to\mathbb{Z}G$ given by $\bar{\lambda}(g)=w(g)\cdot \lambda(g^{-1})$
extends to an antiisomorphism $\bar{ }:\widehat{\mathbb{Z}G}_\chi\to \widehat{\mathbb{Z}G}_{-\chi}$. By
(\ref{orienter}) we now get
\begin{eqnarray}\label{incirel}
[p:q]_v&=&(-1)^i\overline{[q:p]_{-v}}.
\end{eqnarray}
If $A$ is a left $\widehat{\mathbb{Z}G}_{-\chi}$ module, we can turn
${\rm Hom}_{\widehat{\mathbb{Z}G}_{-\chi}}(A,\widehat{\mathbb{Z}G}_{-\chi})$ into a left
$\widehat{\mathbb{Z}G}_\chi$ module by setting $\lambda \cdot \varphi:a\mapsto\varphi(a)\cdot\bar{\lambda}\in
\widehat{\mathbb{Z}G}_{-\chi}$.\\
Let $C^\ast(-\omega,-v)={\rm Hom}_{\widehat{\mathbb{Z}G}_{-\chi}}(C_\ast(-\omega,-v),\widehat{
\mathbb{Z}G}_{-\chi})$. Using (\ref{incirel}) it is easy to see that
\[
\begin{array}{crcl}
P:&C_\ast(\omega,v)&\longrightarrow&C^{n-\ast}(-\omega,-v)\\
 &p&\mapsto&(-1)^{i(i+1)/2}\bar{p}\end{array}\]
is a simple isomorphism of free $\widehat{\mathbb{Z}G}_\chi$ chain complexes, where $\bar{p}:C_\ast(-\omega,
-v)\to \widehat{\mathbb{Z}G}_{-\chi}$ is defined by $\bar{p}(p)=1$ and 0 for all other critical points. This
induces the Poincar\'{e} duality isomorphism $P_i:H_i(M;\widehat{\mathbb{Z}G}_\chi)\to H^{n-i}(M,
\widehat{\mathbb{Z}G}_{-\chi})$.\\
To get a duality isomorphism for the noncommutative localization $\Sigma^{-1}_\chi\mathbb{Z}G$ we need the
following lemma.
\begin{lem}\label{antiinv}
Let $R$ be a ring with unit, $\bar{ }:R\to R$ an involution, $\Sigma$ a set of diagonal matrices over $R$
which is closed under transpose. Then the involution extends to an antiisomorphism $\bar{ }:\Sigma^{-1}R
\to \bar{\Sigma}^{-1}R$.
\end{lem}
\begin{proof}
For any ring $S$ denote $S^o$ the opposite ring, i.e.\ multiplication is given by $(x,y)\mapsto y\cdot x$.
Hence we can think of the involution as a ring homomorphism $\varphi:R^o\to R$. Let $\bar{\varepsilon}:R\to
\bar{\Sigma}^{-1}R$ be the natural map. Now $\bar{\varepsilon}\circ \varphi(\Sigma)\subset GL(\bar{\Sigma}
^{-1}R)$. Note that $A\in\Sigma$ can be thought of as a matrix over $R^o$ and then $\varphi(A)$ is a matrix
over $R$ contained in $\bar{\Sigma}$. Therefore we have a unique map $\theta_1:\Sigma^{-1}R^o\to\bar{\Sigma}
^{-1}R$ such that $\theta_1\circ\varepsilon'=\bar{\varepsilon}\circ\varphi$ with $\varepsilon':R^o\to
\Sigma^{-1}R^o$ the natural map. Similarly we get a unique map $\theta_2:\bar{\Sigma}^{-1}R\to\Sigma^{-1}
R^o$ such that $\varepsilon'\circ\varphi^o=\theta_2\circ\bar{\varepsilon}:R\to\Sigma^{-1}R^o$. It follows
that $\theta_1$ and $\theta_2$ are mutually inverse isomorphisms.\\
We have $\varepsilon'(R^o)\subset\Sigma^{-1}R^o$, so $(\varepsilon')^o(R)\subset(\Sigma^{-1}R^o)^o$. Also
if $A\in\Sigma$, then $A^T\in\Sigma$ and $\varepsilon'(A^T)$ is invertible in $\Sigma^{-1}R^o$. But if a
matrix is invertible over a ring $S$, its transpose is invertible over $S^o$. Therefore $(\varepsilon')^o(A)$
is invertible in $(\Sigma^{-1}R^o)^o$. Thus there is a ring homomorphism $\psi_1:\Sigma^{-1}R\to(\Sigma^{-1}
R^o)^o$ such that $(\varepsilon')^o=\psi_1\circ\varepsilon$ where $\varepsilon:R\to \Sigma^{-1}R$ is the
natural map. Similarly we get a unique ring homomorphism $\psi_2:\Sigma^{-1}R^o\to(\Sigma^{-1}R)^o$ with
$\psi_2\circ\varepsilon'=\varepsilon^o:R^o\to(\Sigma^{-1}R)^o$. It follows that $\psi_1$ and $\psi_2^o$ are
mutually inverse isomorphisms. Now $\theta_1\circ\psi_1^o:(\Sigma^{-1}R)^o\to\bar{\Sigma}^{-1}R$ induces
the desired antiisomorphism.
\end{proof}
Now let $P:C_\ast^\Delta(\tilde{M})\to C_\Delta^{n-\ast}(\tilde{M})$ be a Poincar\'e duality simple
chain homotopy equivalence, e.g.\ induced by an exact Morse form $df$. Let $i:\mathbb{Z}G\to
\Sigma^{-1}_\chi\mathbb{Z}G$ be the inclusion. Then we get a simple chain homotopy equivalence
\[{\rm id}\otimes P:C_\ast^\Delta(M;\Sigma^{-1}_\chi\mathbb{Z}G)\to\Sigma^{-1}_\chi\mathbb{Z}G
\otimes_{\mathbb{Z}G}C_\Delta^{n-\ast}(\tilde{M}).\]
Using Lemma \ref{antiinv} we have an isomorphism $\Theta:\Sigma^{-1}_\chi\mathbb{Z}G
\otimes_{\mathbb{Z}G}C_\Delta^{n-\ast}(\tilde{M})\to C_\Delta^{n-\ast}(M;\Sigma^{-1}_{-\chi}\mathbb{Z}G)$
given by $\Theta(r\otimes\sigma^\ast):s\otimes \tau \mapsto s\cdot\sigma^\ast(\tau)\cdot\bar{r}$. Hence
we get a Poincar\'e duality simple chain homotopy equivalence
\[P_i:C_i^\Delta(M;\Sigma_\chi^{-1}\mathbb{Z}G)\to C_\Delta^{n-i}(M;\Sigma_{-\chi}^{-1}\mathbb{Z}G).\]
Because of Poincar\'e duality we now get the following.
\begin{prop} We have:
\begin{enumerate}
\item $C_\ast^\Delta(M;\widehat{\mathbb{Z}G}_\chi)$ is acyclic if and only if
$C_\ast^\Delta(M;\widehat{\mathbb{Z}G}_{-\chi})$ is acyclic.
\item $C_\ast^\Delta(M;\Sigma^{-1}_\chi\mathbb{Z}G)$ is acyclic if and only if
$C_\ast^\Delta(M;\Sigma^{-1}_{-\chi}\mathbb{Z}G)$ is acyclic.
\end{enumerate}
\end{prop}
In that case we get for the Latour obstructions
\[\tau(M,\alpha)=(-1)^{n-1}\bar{\tau}(M,-\alpha)\]
both in ${\rm Wh}(G;\chi)$ and ${\rm Wh}(G;\Sigma_\chi)$ by Milnor \cite{milnwh}. Notice that the
antiisomorphism $\bar{ }:
\widehat{\mathbb{Z}G}_{-\chi}\to \widehat{\mathbb{Z}G}_\chi$ induces an isomorphism of abelian groups
$\bar{ }:{\rm Wh}(G;-\chi)\to{\rm Wh}(G;\chi)$ by taking the conjugate transpose of a matrix. Similar for
${\rm Wh}(G;\Sigma_\chi)$.
\section{Connections between Novikov homology and controlled connectivity}\label{sec8}
Proposition \ref{no0n} and Proposition \ref{no1n-1} show directly how controlled connectivity properties lead
to the vanishing of certain Novikov homology groups and vice versa, at least in the manifold case. In Section
\ref{sec4} we did not deal with end points as we needed absolute $CC^1$ for the results in Section \ref{sec6}.
But we can refine the results of Section \ref{sec4} slightly by looking at end points.\\
For a control function $f$ of $\alpha$ define
\[\tilde{M}^-_t=f^{-1}((-\infty,t])\mbox{ and }\tilde{M}^+_t=f^{-1}([t,\infty))\]
The analogues
of Propositions \ref{cc0con}-\ref{cc1con} are now
\begin{prop}
Let $\alpha\in H^1(M;\mathbb{R})$. Assume that $\alpha\not=0$ and $n=\dim M\geq 3$. Then the following
are equivalent.
\begin{enumerate}
\item $\alpha$ is $CC^0$ at $-\infty$ (resp. $+\infty$).
\item There is a control function $f$ of $\alpha$ without critical points of index 0, $n$ and with connected
$\tilde{M}^-_t$ (resp. $\tilde{M}^+_t$).
\item There is a control function $f$ of $\alpha$ with connected $\tilde{M}^-_t$ (resp.
$\tilde{M}^+_t$).
\end{enumerate}
\end{prop}
The proof is analogous to the proof of Proposition \ref{cc0con}.
\begin{prop}
Let $\alpha\in H^1(M;\mathbb{R})$. Assume that $\alpha\not=0$ and $n=\dim M\geq 5$. Then $\alpha$ is $CC^0$
at $-\infty$ (resp. $+\infty$) if and only if $\alpha$ can be represented by a Morse form $\omega$ without
critical points of index 0, 1 and $n$ (resp. 0, $n-1$ and $n$).
\end{prop}
\begin{proof}
Replace $\tilde{N}(f,t)$ by $\tilde{M}^-_t$ in the proof of Proposition \ref{no1n-1}, the rest is analogous.
\end{proof}
\begin{prop}
Let $\alpha\in H^1(M;\mathbb{R})$. Assume that $\alpha\not=0$ and $n=\dim M\geq 5$. Then the following are
equivalent.
\begin{enumerate}
\item $\alpha$ is $CC^1$ at $-\infty$ (resp. $+\infty$).
\item There is a control function $f$ of $\alpha$ without critical points of index 0, 1 and $n$ (resp.
0, $n-1$ and $n$) and with simply connected $\tilde{M}^-_t$ (resp. $\tilde{M}^+_t$).
\item There is a control function $f$ of $\alpha$ with simply connected $\tilde{M}^-_t$
(resp. $\tilde{M}^+_t$).
\end{enumerate}
\end{prop}
\begin{exam}\em \label{bausol}
Let $M$ be a closed connected smooth manifold such that its fundamental group is the Baumslag-Solitar
group $G=\langle x,t\,|\,t^{-1}xt=x^2\rangle$. Clearly $H_1(M)=\mathbb{Z}$. Let $\alpha\in H^1(M;\mathbb{R})$
induce the homomorphism $\chi:G\to\mathbb{Z}$ given by $x\mapsto 0$ and $t\mapsto 1$. It is shown in
\cite[\S 10.2]{biegeo} that $\chi$ is $CC^1$ at $-\infty$, but not $CC^0$ at $+\infty$.
\end{exam}
This shows that we can find cohomology classes which are $CC^1$ over $-\infty$ but not
$CC^0$ over $+\infty$. In particular we can represent such a cohomology class by a Morse form without
critical points of index 0, 1 and $n$, but with critical points of index $n-1$.\\
Let us now return to the group theoretic setting. Given a character $\chi:G\to\mathbb{R}$ let $X$ again be
the $k$-skeleton of the universal cover of a $K(G,1)$ CW-complex with finite $k$-skeleton and $h$ a
control function. We can look at the completed cellular complex $\widehat{\mathbb{Z}G}_\chi
\otimes_{\mathbb{Z}G}C_\ast(X)$ and the completed singular complex $\widehat{\mathbb{Z}G}_\chi
\otimes_{\mathbb{Z}G}C_\ast^s(X)$ and denote its homology by $H_\ast(X;\widehat{\mathbb{Z}G}_\chi)$. For
this situation let us introduce a notion similar to controlled connectivity.
\begin{defi}\em
The homomorphism $\chi:G\to\mathbb{R}$ is called \em controlled $(k-1)$-acyclic ($CA^{k-1}$) over $-\infty$,
\em if for every $s\in\mathbb{R}$ and $p\leq k-1$ there is an $\lambda(s)\geq 0$ such that every singular
$p$-cycle (over $\mathbb{Z}$) in $X_s$ bounds in $X_{s+\lambda(s)}$ and $s+\lambda(s)\to -\infty$ as $s\to
-\infty$.
\end{defi}
We can define $\chi$ being $CC^{k-1}$ over $+\infty$ similarly. For $k\leq 1$ we clearly have $\chi$ is
$CC^{k-1}$ over $-\infty$ if and only if $\chi$ is $CA^{k-1}$ over $-\infty$. For higher $k$ we have the
usual problem in comparing homology and homotopy, but there is a Hurewicz-type theorem, see Geoghegan
\cite{geoghe}.
\begin{theo}
For $k\geq 2$, $\chi$ is $CC^{k-1}$ over $-\infty$ if and only if $\chi$ is $CC^1$ over $-\infty$ and
$CA^{k-1}$ over $-\infty$.
\end{theo}
The relation with Novikov homology is now summarized in
\begin{prop}{\rm \cite[Prop.D.2]{bieri}}\qua
\label{propol}
Let $\chi:G\to\mathbb{R}$ be a character, $k\geq 1$ and $X$ as
above. Then $\chi$ is $CA^{k-2}$ over $-\infty$ if and only if
$H_i(X;\widehat{\mathbb{Z}G}_\chi)=0$ for $i\leq k-1$.
\end{prop}
\begin{proof}
We can attach $(k+1)$-cells to $X$ to make $X$ $k$-connected. This will not change the Novikov homology in
dimensions $\leq k-1$.
We can describe $CA^{k-1}$ by saying that the map $\tilde{H}_i(X_s)\to \tilde{H}_i(X_{s+\lambda(s)})$ induced
by inclusion is trivial for $i\leq k-1$ with $s$ and $\lambda(s)$ as in the definition. We have the
commutative diagram
\[\begin{array}{ccc}
H_{i+1}(X,X_s)&\longrightarrow&\tilde{H}_i(X_s)\\[0.3cm]
\Big\downarrow& &\big\downarrow\\[0.3cm]
H_{i+1}(X,X_{s+\lambda(s)})&\longrightarrow&\tilde{H}_i(X_{s+\lambda(s)})
\end{array}
\]
and the horizontal arrows are isomorphisms for $i\leq k-1$ since $X$ is $k$-connected.\\
It is known that $\widehat{\mathbb{Z}G}_\chi\otimes_{\mathbb{Z}G}C^s_\ast(X)=\lim\limits_{\longleftarrow}
C_\ast^s(X,X_s)$, compare Remark \ref{remarkos}, so the Novikov homology fits into a short exact sequence
\[0\longrightarrow{\lim\limits_{\longleftarrow}}^1 H_{i+1}(X,X_s)\longrightarrow H_i(X,
\widehat{\mathbb{Z}G}_\chi)\longrightarrow\lim\limits_{\longleftarrow} H_i(X,X_s)\longrightarrow 0\]
see e.g.\ Geoghegan \cite{geoghe}. By the diagram above and this short
exact sequence we now get immediately that $CA^{k-1}$ implies the
vanishing of the Novikov homology groups in dimensions $\leq k-1$ and
this vanishing implies $CA^{k-2}$. To see that already $CA^{k-2}$
implies $H_{k-1}(X;\widehat{\mathbb{Z}G}_\chi)=0$ note that by Bieri
and Renz \cite[Th.4.2]{bieren} the inverse system $\{H_k(X,X_s)\}$ is
surjective, hence ${\lim\limits_{\longleftarrow}}^1 H_k(X,X_s)=0$. By
the short exact sequence above we get the result.
\end{proof}
Let us now look at the case of an aspherical manifold $M$. In this case we can use the universal cover
$\tilde{M}$ to check for all controlled connectivity properties.
\begin{prop}
Let $M$ be an aspherical closed connected smooth manifold with $n=\dim M$ and $\chi:G\to\mathbb{R}$ a
character. Then the following are equivalent.
\begin{enumerate}
\item The Novikov complex
$C_\ast^\Delta(M;\widehat{\mathbb{Z}G}_\chi)$ is acyclic.
\item $\chi$ is $CA^{n-2}$ over $-\infty$.
\item $\chi$ is $CA^{[\frac{n}{2}]-1}$.
\end{enumerate}
\end{prop}
\begin{proof}
By Proposition \ref{propol} we get (1) $\Rightarrow$ (2) and
(1) $\Rightarrow$ (3).\\
If $\chi$ is $CA^{n-2}$ over $-\infty$, we get
$H_i(M;\widehat{\mathbb{Z}G}_\chi)=0$ for $i\leq n-1$ by Proposition
\ref{propol}. Now
$H_n(M;\widehat{\mathbb{Z}G}_\chi)=H^0(M;\widehat{\mathbb{Z}G}_{-\chi})=0$
by Poincar\'e duality and since $\chi\not=0$.\\
If $\chi$ is $CA^{[\frac{n}{2}]-1}$, we get
$H_i(M;\widehat{\mathbb{Z}G}_\chi)=0$ for $i\leq [\frac{n}{2}]$. Now
for $i\geq [\frac{n}{2}]+1$ we have 
\[H_i(M;\widehat{\mathbb{Z}G}_\chi)= H^{n-i}(M;\widehat{\mathbb{Z}G}_{-\chi})\]
But $n-i\leq n- [\frac{n}{2}]-1\leq [\frac{n}{2}]$ and
$H_{n-i}(M;\widehat{\mathbb{Z}G}_{-\chi}) =0$, since $-\chi$ is
$CA^{[\frac{n}{2}]-1}$ as well. Therefore we get the result.
\end{proof}
The proof shows we can loosen the condition that $M$ be aspherical
slightly to get the following.
\begin{coro}
Let $M$ be a closed connected smooth manifold with $n=\dim M$ and
$\chi:G\to \mathbb{R}$ a character such that $\tilde{M}$ is
$[\frac{n}{2}]$-connected. Then
$C_\ast^\Delta(M;\widehat{\mathbb{Z}G}_\chi)$ is acyclic if and only if
$\chi$ is $CA^{[\frac{n}{2}]-1}$.
\end{coro}
For an aspherical manifold $M$ Latour's theorem can now be phrased as follows.
\begin{theo}\label{alat}
Let $M$ be an aspherical closed connected smooth manifold with $n=\dim M\geq 6$ and $\chi:G\to\mathbb{R}$
a character. Then $\chi$ can be represented by a nonsingular closed 1-form if and only if $\chi$ is $CC^1$,
$\chi$ is $CA^{n-2}$ over $-\infty$ and $\tau(M,\chi)=0$.
\end{theo}
Whitehead groups of aspherical manifolds are conjectured to be zero which is known for certain classes of
manifolds. In this case $CC^1$ and $CA^n$ over $-\infty$ suffices in Theorem \ref{alat}.
\begin{appendix}
\section{Chain homotopy equivalences between Novikov complexes}\label{apchain}
In this appendix we introduce several chain homotopy equivalences between Novikov complexes and sketch
proofs of their properties. The techniques involved are described in more detail in \cite[App.A]{schuet}
and \cite[\S 9]{schue2}. The reader might also want to compare Cornea and Ranicki \cite{corran}, Hutchings
and Lee \cite[\S 2.3]{hutlee},
Latour \cite[\S 2]{latour}, Pozniak \cite[\S 2]{poznia} and Schwarz \cite{schwab,schwar}.
\subsection*{The Morse-Smale complex}
Let us begin with the exact case. Let $(W;M_0,M_1)$ be a compact cobordism, $f:W\to\mathbb{R}$ a Morse
function and $v$ an $f$-gradient satisfying the transversality condition. A smooth triangulation $\Delta$ of
$W$ is said to be \em adjusted to \em $v$, if every $k$-simplex $\sigma$ intersects the unstable manifolds
$W^u(p,v)$ transversely for all critical points $p$ of index $\geq k$. In particular, if $p$ is a critical
point of index $k$, a $k$-simplex $\sigma$ intersects $W^u(p,v)$ in finitely many points. Using the
orientations we can assign to every such point a sign. Given a regular covering space $q:\tilde{W}\to W$ we
can use the covering transformation group $G$ and liftings of critical points and simplices to
assign an element $[\sigma:p]\in\mathbb{Z}G$ to the intersection and define a map: \\
\begin{minipage}{3.5cm}\begin{eqnarray*}\label{chheq}\\ \nonumber\end{eqnarray*}
\end{minipage}\begin{minipage}{8cm}
\begin{eqnarray*}
\varphi_v:C^\Delta_\ast(\tilde{W},\tilde{M}_0)&\longrightarrow& C^{MS}_\ast(\tilde{W},\tilde{M}_0,f,v)\\
\sigma_k&\mapsto&\sum\limits_{p\in{\rm crit}_k(f)}[\sigma:p]\,p  \\
\end{eqnarray*}\end{minipage}

Here $C^{MS}_\ast(\tilde{W},\tilde{M}_0,f,v)$ is the Morse-Smale complex generated by the critical points
of $f$. For $A\subset W$ we denote $\tilde{A}=q^{-1}(A)$. It is shown in \cite[App.A]{schuet} that
adjusted triangulations are generic and $\varphi_v$ is a simple homotopy equivalence.\\
Now given another Morse function $g:W\to\mathbb{R}$ with a transverse $g$-gradient $w$, let
$\Phi:W\to W$ be isotopic to the identity such that $\Phi(W^s(q,v))\,\pitchfork W^u(p,w)$ for
critical points $q$ of $f$ and $p$ of $g$ with ${\rm ind}\,q\leq{\rm ind}\,p$. The existence of $\Phi$ is
achieved by standard transversality arguments. Furthermore we get openness and density for such $\Phi$ in
the smooth topology. If ${\rm ind}\,q={\rm ind}\,p$ we get that $\Phi(W^s(q,v))\cap W^u(p,w)$ is finite, in
fact we get an intersection number $[q:p]\in\mathbb{Z}G$ as above and we can define $\psi_{v,w}:
C_\ast^{MS}(\tilde{W},\tilde{M}_0,f,v)\to C_\ast^{MS}(\tilde{W},\tilde{M}_0,g,w)$ by
\[\psi_{v,w}(q)=\sum_{p,{\rm ind}\,p={\rm ind}\,q}[q:p]\,p.\]
The proof that $\psi_{v,w}$ is a chain map is identical to \cite[\S 9]{schue2}, even though the two Morse
functions there were equal. Also, as in \cite[Lm.A.2]{schuet} the chain homotopy type does not depend on
$\Phi$.
\begin{prop}\label{exaprop}
For $i=0,1,2$ let $f_i:W\to\mathbb{R}$ be a Morse function of the cobordism $(W;M_0,M_1)$ and $v_i$ a
transverse $f_i$-gradient. Then
\begin{enumerate}
\item $\psi_{v_0,v_1}\circ\varphi_{v_0}\simeq\varphi_{v_1}$
\item $\psi_{v_1,v_2}\circ\psi_{v_0,v_1}\simeq\psi_{v_0,v_2}$.
\end{enumerate}
\end{prop}
In particular we get that $\psi_{v,w}$ is a simple chain homotopy equivalence.
\begin{proof}
The proof of (1) is identical to the proof of \cite[Prop.9.4]{schue2} even though the Morse functions there
are equal. (2) now follows from the fact that $\varphi_{v_i}$ is a chain homotopy equivalence, but in view
of the nonexact case let us give a direct proof. Let $\Phi:W\to W$ be isotopic to the identity such that 
\begin{equation}\label{2chheq}
\begin{array}{rcl}
\Phi(W^s(q,v_0))&\pitchfork & W^u(p,v_1)\\
\Phi(W^s(q,v_0))&\pitchfork & W^u(r,v_2)\\
\Phi(W^s(p,v_1))&\pitchfork & W^u(r,v_2)\end{array}
\end{equation}
for the relevant critical points. For $j=-1,\ldots,n$ and $\delta>0$ let
\[D^j_\delta(v_i)=\bigcup_{p\in{\rm crit}\,f_i \atop {\rm ind}\,p\leq j} D_\delta(p,v_i)\cup M_0.\]
Choose $\delta>0$ so small that $\Phi(D_\delta^j(v_i))$ is disjoint from $W^u(p,v_k)$ where $k>i$ and
${\rm ind}\,p>j$. This is possible by (\ref{2chheq}).\\
Let $\Theta:W\times\mathbb{R}\to W$ be induced by the flow of $-v_1$, i.e.\ stop once the boundary is
reached. There is a $K>0$ such that $\Phi(\Theta_K(D^j_\delta(v_0))\subset D^j_\delta(v_1)$. Let $h:W\times
I\to W$ be a homotopy between the identity and $\Phi\circ\Theta_K$ such that $\Phi(h(W^s(p,v_0)\times I))
\,\pitchfork W^u(r,v_2)$ for ${\rm ind}\,p\leq{\rm ind}\,r-1$. Again we get intersection numbers
$[p:r]\in\mathbb{Z}G$. Then $h$ defines a chain homotopy $H:C^{MS}_\ast(\tilde{W},\tilde{M}_0,f_0,v_0)\to
C_{\ast+1}^{MS}(\tilde{W},\tilde{M}_0,f_2,v_2)$ between $\psi_{v_1,v_2}\circ\psi_{v_0,v_1}$ and
$\psi_{v_0,v_2}$ by
\[H(p)= (-1)^{{\rm ind}\,p}\sum_{r,\,{\rm ind}\,r={\rm ind}\,p+1}[p:r]\,r.\]
To see that this is indeed the right chain homotopy compare the proof of \cite[Prop.9.4]{schue2}.
\end{proof}
\subsection*{The Novikov complex}
Let $M$ be a closed connected smooth manifold and $\omega_i$ be cohomologous Morse forms with transverse
$\omega_i$-gradients $v_i$ for $i=0,1$. Then we can define chain maps
\[\varphi_{v_i}:C^\Delta_\ast(M;\widehat{\mathbb{Z}G}_\chi)\to C_\ast(\omega_i,
v_i)\]
and
\[\psi_{v_0,v_1}:C_\ast(\omega_0,v_0)\to C_\ast(\omega_1,v_1)\]
as in the exact case using intersection numbers which are now elements of $\widehat{\mathbb{Z}G}_\chi$.
To see this one uses inverse limit arguments in the rational case, compare the proof of
\cite[Prop.9.2]{schue2}. The irrational case is treated by approximation, one shows that $[\sigma:q]$ and
$[p:q]$ are elements of $\widehat{\mathbb{Z}G}_\chi\cap\widehat{\mathbb{Z}G}_{\chi'}$, where
$\chi':G\to\mathbb{Q}$. The details are similar to \cite[Prop.9.2]{schue2}, though the Morse form is fixed
there, and will be omitted.
\begin{prop}\label{formcom}
For $i=0,1,2$ let $\omega_i$ be cohomologous Morse forms and $v_i$ transverse $\omega_i$-gradients. Then
\begin{enumerate}
\item $\psi_{v_0,v_1}\circ\varphi_{v_0}\simeq\varphi_{v_1}$.
\item $\psi_{v_1,v_2}\circ\psi_{v_0,v_1}\simeq\psi_{v_0,v_2}$.
\end{enumerate}
\end{prop}
\begin{proof}
Both statements are deduced from the exact case by inverse limit arguments in the rational and approximation
arguments in the irrational case. Compare the proof of \cite[Prop.9.5]{schue2}.
\end{proof}
\begin{coro}\label{chainequ}
$\psi_{v_0,v_1}$ and $\varphi_{v_0}$ are chain homotopy equivalences.
\end{coro}
\begin{proof}
That $\psi_{v_0,v_1}$ is a chain homotopy equivalence follows from Proposition \ref{formcom}.2 since
$\psi_{v_0,v_0}\simeq{\rm id}$. To see that $\varphi_{v_0}$ is a chain homotopy equivalence, it is by
Proposition \ref{formcom}.1 good enough to find a $v_1$ such that $\varphi_{v_1}$ is a chain homotopy
equivalence. But by a nice trick of Latour \cite[Lm.2.28]{latour} there is a Morse form $\omega_1$
cohomologous to $\omega_0$ and a transverse $\omega_1$-gradient $v_1$ such that $v_1$ is also the gradient
of an ordinary Morse function $f:M\to\mathbb{R}$. Then $C_\ast(\omega_1,v_1)=\widehat{\mathbb{Z}G}_\chi
\otimes_{\mathbb{Z}G}C_\ast(\tilde{M},f,v_1)$ and $\varphi_{v_1}={\rm id}_{\widehat{\mathbb{Z}G}_\chi}
\otimes_{\mathbb{Z}G}\varphi^{MS}_{v_1}$. Since $\varphi^{MS}_{v_1}$ is a chain homotopy equivalence, so is
$\varphi_{v_1}$.
\end{proof}
We are also interested in torsion.
\begin{theo}\label{simple}
$\psi_{v_0,v_1}$ and $\varphi_{v_0}$ are simple chain homotopy equivalences, i.e.\ $\tau(\psi_{v_0,v_1})=
\tau(\varphi_{v_0})=0\in{\rm Wh}(G;\chi)$.
\end{theo}
\begin{proof}
That $\tau(\varphi_{v_0})$ is in the image of units of the form $1-a$ with $\|a\|<1$ is shown in
\cite{schue2}. Now $\tau(\psi_{v_0,v_1})=0$ follows from Proposition \ref{formcom}.1.\\
Alternatively we can use the techniques of Latour \cite[\S 2.25-2.28]{latour} to show that
$\tau(\psi_{v_0,v_1})=0$. Then $\tau(\varphi_{v_0})=0$ follows from Proposition \ref{formcom}.1 after
noticing that $\tau(\varphi_{v_1})=0$ for $\varphi_{v_1}$ as in the proof of Corollary \ref{chainequ}.
\end{proof}
\begin{rema}\em\label{remarkos}
Both proofs that $\tau(\varphi_v)=0$ are quite involved. But in the rational case there is a significantly
easier proof: let $\rho:M_\omega\to M$ be the infinite cyclic covering such that $\rho^\ast\omega=df$.
We can assume that 0 is a regular value of $f$ and that $f(tx)-f(x)=1$ for a generator $t$ of the infinite
cyclic covering transformation group. For $k$ a positive integer let $M_k=f^{-1}([-k,0])$ and $N_k=f^{-1}
(\{-k\})$. Then the following diagram commutes
\[\begin{array}{ccc}
C^\Delta_\ast(\tilde{M}_k,\tilde{N}_k)&\longleftarrow&C_\ast^\Delta(\tilde{M}_{k+1},\tilde{N}_{k+1})\\[0.3cm]
\Big\downarrow\varphi_{v|}& &\big\downarrow\varphi_{v|}\\[0.3cm]
C_\ast^{MS}(\tilde{M}_k,\tilde{N}_k,f|,v|)&\longleftarrow&C_\ast^{MS}(\tilde{M}_{k+1},\tilde{N}_{k+1})
\end{array}\]
Let $\widehat{\mathbb{Z}G}_\chi^0$ be the subring of $\widehat{\mathbb{Z}G}_\chi$ consisting of elements $a$
with $\|a\|\leq 1$ and let $H=\ker\chi$. Then the inverse
limits are finitely generated free $\widehat{\mathbb{Z}G}_\chi^0$ complexes. Since $\varphi_{v|}$ is a
chain homotopy equivalence, so is $\lim\limits_{\longleftarrow}\varphi_{v|}$. Also ${\rm id}_{\mathbb{Z}H}
\otimes_{\widehat{\mathbb{Z}G}_\chi^0}\lim\limits_{\longleftarrow}\varphi_{v|}=\varphi_{v|}:C_\ast^\Delta
(\tilde{M}_1,\tilde{N}_1)\to C_\ast^{MS}(\tilde{M}_1,\tilde{N}_1,f|,v|)$ and ${\rm id}_{\widehat{\mathbb{Z}
G}_\chi}\otimes_{\widehat{\mathbb{Z}G}_\chi^0}\lim\limits_{\longleftarrow}\varphi_{v|}=\varphi_v:
C^\Delta_\ast(M;\widehat{\mathbb{Z}G}_\chi)\to C_\ast(\omega,v)$. Since
$\varphi_{v|}$ is simple, $\tau(\lim\limits_{\longleftarrow}\varphi_{v|})\in\ker r_\ast\subset K_1(
\widehat{\mathbb{Z}G}_\chi^0)/\langle\pm[h]\,|\,h\in H\rangle$, where $r:\widehat{\mathbb{Z}G}_\chi^0\to
\mathbb{Z}H$ is projection. But by an elementary argument $\ker r_\ast$ is generated by units of the form
$1-a$ with $\|a\|<1$, see Pajitnov \cite[Lm.1.1]{pajito}. Hence $\tau(\varphi_v)=i_\ast\tau(\lim
\limits_{\longleftarrow}\varphi_{v|})=0\in{\rm Wh}(G;\chi)$ where $i:\widehat{\mathbb{Z}G}_\chi^0\to
\widehat{\mathbb{Z}G}_\chi$ is inclusion.
\end{rema}
This proof does not seem to carry over to the irrational case.
\subsection*{Continuation}
Given two Morse-Smale or Novikov complexes, one can find other methods in the literature to produce a
chain homotopy equivalence between these complexes, like continuation. This principle is explained
e.g.\ in Schwarz \cite{schwab} or Pozniak \cite[\S 2]{poznia}. The purpose of this subsection is to show
that even though its definition differs from the definition of $\psi_{v,w}$ given above it agrees with
$\psi_{v,w}$ up to chain homotopy. We will only consider the exact case noting that the nonexact case can
be derived from the exact case by the typical techniques described above. To describe continuation we
choose the description of Pozniak \cite[\S 2.6]{poznia}.\\
So let $f,g:M\to\mathbb{R}$ be Morse functions, $v$ a transverse $f$-gradient and $w$ a transverse
$g$-gradient. Let $F:M\times\mathbb{R}\to\mathbb{R}$ be a smooth function with $F(x,0)=f(x)+C_1$, $F(x,1)
=g(x)+C_2$ such that the critical points of $F$ are exactly of the form $(p,0)$ where $p$ is a critical
point of $f$ and $(q,1)$ where $q$ is a critical point of $g$. Furthermore we want ${\rm ind}\,(p,0)=
{\rm ind}\,p+1$ and ${\rm ind}\,(q,1)={\rm ind}\,q$. Using a transverse $F$-gradient $u$ which agrees with
$v$ on $M\times\{0\}$ and with $w$ on $M\times\{1\}$, Pozniak \cite[\S 2.6]{poznia} shows that there is
an acyclic finitely generated free Morse-Smale complex $C_\ast^{MS}(F,u)$ which fits into a short exact
sequence of chain complexes
\[0\longrightarrow C_\ast^{MS}(\tilde{M},g,w)\longrightarrow C_\ast^{MS}(F,u)\longrightarrow
C^{MS}_{\ast-1}(\tilde{M},f,v)\longrightarrow 0\]
But this means we can think of $C_\ast^{MS}(F,u)$ as the mapping cone of a chain homotopy equivalence
$c_{v,w}:C^{MS}_\ast(\tilde{M},f,v)\to C_\ast^{MS}(\tilde{M},g,w)$. Furthermore $c_{v,w}$ can be described
by flowlines of $-u$ from critical points $(p,0)$ to critical points $(q,1)$. Notice that this agrees
with the chain map given in Cornea and Ranicki \cite[Prop.1.11]{corran}.
\begin{prop}
We have $c_{v,w}\simeq \psi_{v,w}: C^{MS}_\ast(\tilde{M},f,v)\to C_\ast^{MS}(\tilde{M},g,w)$.
\end{prop}
\begin{proof}
We can assume that $W^s(p,v)\,\pitchfork W^u(q,w)$ for ${\rm ind}\,p\leq{\rm ind}\,q$. Let $p$ be
a critical point of $f$ of index $i$. Let $\theta_p:\mathbb{R}^i\times\mathbb{R}\to W^s((p,0),u)$ be an
immersion of the stable manifold in $M\times\mathbb{R}$ so that we can identify the image of $\mathbb{R}^i
\times \{0\}$ with $W^s(p,v)$ in $M=M\times \{0\}$. By the definition of $F$ we either have $\theta_p
(x,t)\in M\times (0,1)$ for all $x\in\mathbb{R}^i$ and $t>0$, or for all $x\in\mathbb{R}^i$ and $t<0$. Let us
assume this is true for $t>0$.\\
Identify $C_k^{MS}(\tilde{M},g,w)=H_k(C^k(w),C^{k-1}(w))$ where 
\[C^k(w)=M-\bigcup\limits_{{\rm ind}\,q>k}
W^u(q,w),\] 
compare \cite[\S 9]{schue2}. By the transversality assumption on $u$ we can find for every
compact disc $D^i\subset \mathbb{R}^i$ a $K>0$ such that $p_M\circ\theta_p(D^i\times\{K\})\subset C^i(w)$,
since $\theta_p(\mathbb{R}^{i+1})$ will avoid critical points $(q,1)$ with ${\rm ind}\,q>i$. Here $p_M:
M\times\mathbb{R}\to M$ is projection. We can also find a large disc $D^i_p\subset\mathbb{R}^i$ such
that $p_M\circ\theta_p(\partial D_p^i\times\{K\})\subset C^{i-1}(w)$ and
\[c_{v,w}(p)=(p_{\tilde{M}}\circ\tilde{\theta}_p)_\ast[D^i_p\times\{K\}]\in H_i(\tilde{C}^i(w),
\tilde{C}^{i-1}(w)).\]
If $D^i_p$ is large enough, we also have
\[\psi_{v,w}(p)=(p_{\tilde{M}}\circ\tilde{\theta}_p)_\ast[D^i_p\times \{0\}].\]
Choose $K$ so large that it works for every critical point of $f$. Let $C$ be the union of the images of the
discs $D^i_p\times\{0\}$ in $M$ and let $h:C\times I\to M$ be a homotopy between $p_M\circ\theta(D^i_p
\times\{0\})$ and $p_M\circ\theta_p(D^i_p\times\{K\})$ such that $h(p_M\circ\theta_p(D^i_p\times \{0\})
\times I)$ intersects $W^u(q,w)$ transversely for ${\rm ind}\,p\leq{\rm ind}\,q-1$. Then $h$ extends to a
homotopy $h:M\times I\to M$ of the identity which induces the desired chain homotopy equivalence.
\end{proof}
\end{appendix}

\Addressesr

\begin{thebibliography}
\bibitem{bieri}R. Bieri, \textsl{The geometric invariants of a group: a survey
with emphasis on the homotopical approach}, Geometric group theory, Vol. 1
(Sussex, 1991), 24-36, London Math. Soc. Lecture Note Ser. 181, Cambridge Univ. Press, Cambridge, 1993.
\bibitem{biegeo}R. Bieri and R. Geoghegan, \textsl{Connectivity properties of group actions on non-positively curved
spaces}, preprint, to appear in Mem. Amer. Math. Soc.
\bibitem{biege2}R. Bieri and R. Geoghegan, \textsl{Kernels of actions on
non-positively curved spaces}, Geometry and cohomology in group theory (Durham,
1994), 24-38, London Math. Soc. Lecture Note Ser. 252, Cambridge Univ. 
Press, Cambridge, 1998.
\bibitem{binest}R. Bieri, W. Neumann and R. Strebel, \textsl{A geometric invariant of discrete groups}, Invent.
Math. 90 (1987), 451-477.
\bibitem{bieren}R. Bieri and B. Renz, \textsl{Valuations on free resolutions and higher geometric invariants of
groups}, Comment. Math. Helv. 63 (1988), 464-497.
\bibitem{brown}K. Brown, \textsl{Finiteness properties of groups}, J. Pure Appl. Algebra 44 (1987), 45-75.
\bibitem{cohn}P. Cohn, \textsl{Free rings and their relations}. Second edition. London Mathematical Society
Monographs, 19. Academic Press, Inc., London, 1985.
\bibitem{corran}O. Cornea and A. Ranicki, \textsl{Rigidity and glueing for the Morse and Novikov complexes}, preprint,
available as math.AT/0107221.
\bibitem{damian}M. Damian, \textsl{Formes ferm\'ees non singuli\`eres et propri\'et\'es de finitude des groupes},
Ann. Sci. \'Ecole Norm. Sup. (4) 33 (2000), 301-320.
\bibitem{dold}A. Dold, \textsl{Lectures on algebraic topology}, Die Grundlehren der mathematischen Wissenschaften, Band 200. Springer-Verlag, New
York-Berlin, 1972.
\bibitem{farbef}M. Farber, \textsl{Sharpness of the Novikov inequalities}, Funktsional. Anal. i Prilozhen. 19 (1985),
49-59. English translation in Functional Anal. Appl. 19 (1985), 40-48.
\bibitem{farber}M. Farber, \textsl{Morse-Novikov critical point theory, Cohn localization and Dirichlet units},
Commun. Contemp. Math. 1 (1999), 467-495.
\bibitem{farran}M. Farber and A. Ranicki, \textsl{The Morse-Novikov theory of circle-valued functions and
noncommutative localization}, Tr. Mat. Inst. Steklova 225 (1999) 381-388.
\bibitem{farrel}F. Farrell, \textsl{The obstruction to fibering a manifold over a circle}, Yale University Ph.D.
thesis, 1967.
\bibitem{farre2}F. Farrell, \textsl{The obstruction to fibering a manifold over a circle},
Actes du Congr\`es International des Math\'ematiciens (Nice, 1970), Tome 2, 69-72, Gauthier-Villars, 1971.
\bibitem{geoghe}R. Geoghegan, \textsl{Topological methods in group theory}, monograph in preparation.
\bibitem{hutlee}M. Hutchings and Y-J. Lee, \textsl{Circle-valued Morse theory, Reidemeister torsion, and
Seiberg-Witten invariants of three manifolds}, Topology 38 (1999), 861-888.
\bibitem{latour}F. Latour, \textsl{Existence de 1-formes ferm\'ees non singuli\`eres dans une classe de cohomologie
de de Rham}, Publ. IHES No.80 (1994), 135-194.
\bibitem{milnhc}J. Milnor, \textsl{Lectures on the h-cobordism theorem}, 
 Notes by L. Siebenmann and J. Sondow, Princeton University Press, Princeton, 
N.J., 1965.
\bibitem{milnwh}J. Milnor, \textsl{Whitehead torsion}, Bull. Amer. Math. Soc. 72 (1966), 358-426.
\bibitem{noviko}S. Novikov, \textsl{Multivalued functions and functionals. An analogue of the Morse theory},
Dokl. Akad. Nauk SSSR 260 (1981), 31-35. English translation in Soviet Math. Dokl. 24 (1981), 222-226.
\bibitem{pajiol}A. Pazhitnov, \textsl{On the sharpness of inequalities of Novikov type for manifolds with a free
abelian fundamental group}, Mat. Sb. 180 (1989), 1486-1523. English translation in Math. USSR-Sb. 68 (1991),
351-389.
\bibitem{pajito}A. Pazhitnov, \textsl{On the Novikov complex for rational Morse forms}, Ann. Fac. Sci. Toulouse 4
(1995), 297-338.
\bibitem{pajisu}A. Pazhitnov, \textsl{Surgery on the Novikov complex}, K-theory 10 (1996), 323-412.
\bibitem{pajisp}A. Pajitnov, \textsl{Incidence coefficients in the Novikov complex for Morse forms: rationality and
exponential growth properties}, available as math.DG/9604004.
\bibitem{pajiov}A. Pajitnov, \textsl{$C^0$-generic properties of boundary operators in the Novikov complex},
Pseudoperiodic topology, Amer. Math. Soc. Transl. Ser. 2, 197 (1999), 29-115.
\bibitem{poznia}M. Pozniak, \textsl{Floer homology, Novikov rings and clean intersections}. Northern California
Symplectic Geometry Seminar, Amer. Math. Soc. Transl. Ser. 2, 196, 1999, 119-181.
\bibitem{ranict}A. Ranicki, \textsl{Finite domination and Novikov rings}, Topology 34 (1995), 619-632.
\bibitem{ranick}A. Ranicki, \textsl{The algebraic construction of the Novikov complex of a circle-valued Morse
function}, available as math.DG/9903090, to appear in Math. Annalen.
\bibitem{renz}B. Renz, Thesis, University of Frankfurt 1987.
\bibitem{schuet}D. Sch\"utz, \textsl{Gradient flows of closed 1-forms and their closed orbits}, available as
math.DG/0009055, to appear in Forum Math.
\bibitem{schue2}D. Sch\"utz, \textsl{One parameter fixed point theory and gradient flows of closed 1-forms}, available
as math.DG/0104245, to appear in K-theory.
\bibitem{schwab}M. Schwarz, \textsl{Morse homology}. Progress in Mathematics, 111. Birkh\"auser Verlag, Basel, 1993.
\bibitem{schwar}M. Schwarz, \textsl{Equivalences for Morse homology}, Geometry and topology in dynamics
(Winston-Salem, NC, 1998/San Antonio, TX, 1999), Contemp. Math. 246 (1999), 197-216.
\bibitem{sharko}V. Sharko, \textsl{The stable algebra of Morse theory}, Izv. Akad. Nauk SSSR Ser. Mat. 54 (1990),
607-631. English translation in Math. USSR-Izv. 36 (1991), 629-653.
\bibitem{sieben}L. Siebenmann, \textsl{A total Whitehead torsion obstruction to fibering over the circle},
Comment. Math. Helv. 45 (1970), 1-48.
\end{thebibliography}
\end{document}